\documentclass[11pt,a4paper]{article}
\usepackage{amssymb,amsmath}
\usepackage{graphicx}
\usepackage{color}
\usepackage{url}
\usepackage{a4wide}

\newcommand {\proof}{\par\noindent{\it Proof}. \ignorespaces}
\def\eproof{\space
        {\ \vbox{\hrule\hbox{\vrule height1.3ex\hskip0.8ex\vrule}\hrule}}
        \bigskip}
\newcommand {\eq} [1] {\begin{equation}\label{#1}}
\newcommand {\en} {\end{equation}}
\newcommand {\eqn}      {\begin{eqnarray}}
\newcommand {\enn}      {\end{eqnarray}}
\newcommand{\bA}        {{\bf A}}
\newcommand{\bB}        {{\bf B}}
\newcommand{\bC}        {{\bf C}}
\newcommand{\bD}        {{\bf D}}
\newcommand{\bE}        {{\bf E}}
\newcommand{\bF}        {{\bf F}}
\newcommand{\bG}        {{\bf G}}
\newcommand{\bH}        {{\bf H}}
\newcommand{\bI}        {{\bf I}}
\newcommand{\bJ}        {{\bf J}}
\newcommand{\bK}        {{\bf K}}
\newcommand{\bL}        {{\bf L}}
\newcommand{\bM}        {{\bf M}}
\newcommand{\bN}        {{\bf N}}
\newcommand{\bP}        {{\bf P}}
\newcommand{\bQ}        {{\bf Q}}
\newcommand{\bR}        {{\bf R}}
\newcommand{\bS}        {{\bf S}}
\newcommand{\bT}        {{\bf T}}
\newcommand{\bU}        {{\bf U}}
\newcommand{\bV}        {{\bf V}}
\newcommand{\bW}        {{\bf W}}

\newcommand{\bY}        {{\bf Y}}
\newcommand{\bZ}        {{\bf Z}}

\newcommand {\bu}       {{\bf u}}

\newcommand {\bw}       {{\bf w}}
\newcommand {\bx}       {{\bf x}}
\newcommand {\by}       {{\bf y}}
\newcommand {\bp}       {{\bf p}}
\newcommand {\bq}       {{\bf q}}
\newcommand {\bz}       {{\bf z}}

\newcommand {\be}       {{\bf e}}
\newcommand{\bfxi}{{\boldsymbol{\boldsymbol{\xi}}}}
\newcommand{\bfeta}{{\boldsymbol{\eta}}}
\newcommand{\bfzeta}{{\boldsymbol{\zeta}}}
\newcommand{\bfphi}{{\boldsymbol{\varphi}}}

\newcommand {\bRn}    {\mathbb{R}^n}
\newcommand {\bRnn}   {\mathbb{R}^{n \times n}}
\newcommand {\bRmn}   {\mathbb{R}^{m \times n}}
\newcommand {\bRnm}   {\mathbb{R}^{n \times m}}
\newcommand {\bRmm}   {\mathbb{R}^{m \times m}}

\newcommand {\mat}      [1] {\left[\begin{array}{#1}}
\newcommand {\rix}          {\end{array}\right]}

\newtheorem{theorem}{Theorem}
\newtheorem{corollary}[theorem]{Corollary}%
\newtheorem{lemma}[theorem]{Lemma}%
\newtheorem{example}[theorem]{Example}%
\newtheorem{remark}[theorem]{Remark}%

\newcommand {\rank}       {\mathop{\rm rank}\nolimits}
\newcommand {\diag}    {\mathop{\rm diag}\nolimits}
\newcommand {\wt}       {\widetilde}
\newcommand {\wh}       {\widehat}
\usepackage[T1]{fontenc}
\usepackage[english]{babel}

\usepackage[inline]{enumitem}
\usepackage{cleveref}
\usepackage{aligned-overset}
\usepackage{xparse}
\usepackage{caption}
\usepackage{csquotes}
\title{Port-Hamiltonian Realizations
 of \\ Nonminimal Linear Time Invariant Systems
}
\author{Christopher Beattie\footnotemark[1], \,  Volker Mehrmann\footnotemark[2], \, and Hongguo Xu\footnotemark[3]
}

\begin{document}
\maketitle

\begin{abstract}
Numerical methods for developing port-Hamiltonian representations of general linear time-invariant systems are studied. The approach extends previous port-Hamiltonian characterizations to include the general non-minimal case and the case where the feedthrough term fails to have an invertible symmetric part.  The resulting construction is able to identify infeasibility when the system fails to be port-Hamiltonian, and allows for the incorporation of perturbations in order to arrive at a nearby port-Hamiltonian system. Results are illustrated via numerical examples.
\end{abstract}
\noindent
{\bf Keywords:}
port-Hamiltonian system, passivity, stability, system transformation, Kalman-Yacubovich-Popov inequality,
Lyapunov inequality, even pencil, quadratic eigenvalue problem.

\noindent
{\bf AMS subject classification.} 93A30, 93B17, 93B11.

\maketitle
\renewcommand{\thefootnote}{\fnsymbol{footnote}}

\footnotetext[1]{
Department of Mathematics, Virginia Tech, Blacksburg, VA 24061, USA.
\texttt{beattie@vt.edu}. Supported by {\it Einstein Foundation Berlin},
through an Einstein Visiting Fellowship.}
\renewcommand{\thefootnote}{\arabic{footnote}}

\footnotetext[2]{
Institut f\"ur Mathematik MA 4-5, TU Berlin, Str. des 17. Juni 136,
D-10623 Berlin, FRG.
\texttt{mehrmann@math.tu-berlin.de}. The authors gratefully acknowledge the support by the Deutsche
Forschungsgemeinschaft (DFG) as part of the collaborative research center SFB 1029 Substantial efficiency increase in gas turbines through direct use of coupled unsteady combustion
and flow dynamics, project A02 {\it Development of a reduced order model of pulsed detonation
combuster} and by Einstein Center ECMath in Berlin.}
\footnotetext[3]{
Department of Mathematics, University of Kansas, Lawrence, KS 66045, USA.
\texttt{xu@math.ku.edu}.
Partially supported by {\it Alexander von Humboldt Foundation}
and by {\it Deutsche Forschungsgemeinschaft},
through the DFG Research Center {\sc Matheon}
{\it Mathematics for Key Technologies}}

\renewcommand{\thefootnote}{\arabic{footnote}}
\setcounter{footnote}{0}

\section{Introduction}\label{sec:intro}
{ In this paper we discuss numerical methods for the construction of transformations that bring standard realizations of linear time-invariant systems into port-Hamiltonian form when this is possible and identification of circumstances when this is not possible.
\subsection{Port-Hamiltonian systems}\label{sec:phsystems}}
The synthesis of system models that describe complex physical phenomena often involves the coupling of
independently developed subsystems originating within different disciplines.
Systematic approaches to coupling such diversely generated subsystems prudently follows
a system-theoretic network paradigm that focuses on the transfer of energy, mass, and
other conserved quantities among the subsystems.
When the subsystem models themselves arise from variational principles,  then the aggregate system typically
has structural features that reflects underlying conservation laws and very often it may be characterized
as a \emph{port-Hamiltonian (pH) system}, see \cite{BeaMV19,BeaMXZ18,GolSBM03,MasSB92,MehM19,MehS23,OrtSMM01,OrtSME02,Sch04,Sch06,SchM95,SchM02,SchM13,SchM23} for some major references and \cite{MehU23} for a detailed survey covering also the case of descriptor systems.
Although pH systems may be formulated within a more general framework, we  restrict ourselves to {\em input-state-output pH systems}, which have the  form
 \begin{align} \label{pHdef}
\dot \bx&=\left(\bJ-\bR\right)\nabla_{\!\bx}{\mathcal H}(\bx)+(\bF-\bP)\bu(t), \\
\by(t)&=\ (\bF+\bP)^\top \nabla_{\!\bx}{\mathcal H}(\bx)  + (\bS+\bN) \bu(t), \nonumber
\end{align}
where $\bx: [0,\infty] \to \bRn$ is the $n$-dimensional {\em state vector}; ${\mathcal H}:\bRn\rightarrow [0,\infty)$ is the \emph{Hamiltonian}, a continuously differentiable scalar-valued vector function,
describing the distribution of internal energy among the energy storage elements of the system; $\bJ=-\bJ^\top \in \bRnn$ is the \emph{structure matrix} describing the energy flux among energy storage elements within the system; $\bR=\bR^\top\in \bRnn$ is the \emph{dissipation matrix} describing energy dissipation/loss in the system;  $\bF\pm\bP\in\mathbb{R}^{n\times m}$ are \emph{port} matrices,  describing the manner in which energy enters and exits the system, and $\bS+\bN$, with $\bS=\bS^\top\in \mathbb R^{m\times m}$ and $\bN=-\bN^\top \in \mathbb R^{m\times m}$, describing the direct feed-through of input to output.
The matrices, $\bR$, $\bP$, and $\bS$ must satisfy
\begin{equation} \label{Kdef}
\bK=\left[\begin{array}{lc}
\bR & \bP \\[1mm]
\bP^\top & \bS
\end{array}\right] \geq 0;
\end{equation}
that is, $\bK$ is symmetric positive-semidefinite. This implies, in particular, that $\bR $ and $\bS$ are also positive semidefinite, $\bR\geq 0$ and $\bS\geq 0$.

Port-Hamiltonian systems generalize the classical notion of \emph{Hamiltonian systems} expressed
in our notation as $\displaystyle \dot{\bx}=\bJ\nabla_{\!\bx}{\mathcal H}(\bx)$.
The analog of the \emph{conservation of energy}  for Hamiltonian
systems is for pH systems (\ref{pHdef}), the {\em dissipation inequality}:
\begin{equation}  \label{DissipIneq}
 {\mathcal H}(\bx(t_1))-{\mathcal H}(\bx(t_0)) \leq \int_{t_{0}}^{t_{1}} \by(t)^\top\bu(t)\ dt,
\end{equation}
which has a natural interpretation as asserting that the increase in internal energy of the system, as measured by ${\mathcal H}$, cannot exceed the \emph{total work} done on the system.   ${\mathcal H}(\bx)$ is
a {\em storage function} associated with the {\em supply rate}, $\by(t)^\top\bu(t)$.
In the language of system theory, (\ref{DissipIneq}) constitutes the property that (\ref{pHdef}) is a {\em passive} system \cite{ByrIW91}.

One may verify with elementary manipulations that the  inequality in (\ref{DissipIneq}) is an immediate consequence of the inequality in (\ref{Kdef}), and holds even when the coefficient matrices $\bJ$, $\bR$, $\bF$, $\bP$, $\bS$, and $\bN$ depend on $\bx$ or explicitly on time $t$ (see, \cite{MasSB92}) or, indeed
(with care taken to define suitable operator domains), when they represent linear operators acting on infinite dimensional spaces \cite{JacZ12,SchM02}.   Notice that with a null input, $\bu(t)=0$, the dissipation inequality asserts that ${\mathcal H}(\bx)$ is non-increasing along any unforced system trajectory.  Thus, ${\mathcal H}(\bx)$ defines a Lyapunov function for the unforced system, so pH systems are implicitly Lyapunov stable \cite{HinP05}.  Similarly, ${\mathcal H}(\bx)$ is non-increasing along any system trajectory that produces a null output, $\by(t)=0$, so pH systems also have Lyapunov stable \emph{zero dynamics} \cite{ByrI03}.

Port-Hamiltonian systems constitute a class of systems that is closed under \emph{power-conserving interconnection}.  This means that port-connected pH systems produce an aggregate system that must also be pH.  This aggregate system will then be guaranteed to be both stable and passive.  Modeling with pH systems, thus, represents physical properties in such a way as to facilitate automated modeling \cite{Kle13} while encoding physical properties explicitly into the structure of the equations.  This framework also provides a compelling motivation to identify and preserve pH structure whenever it is present in order to produce high quality reduced order surrogate models, see \cite{BeaG11,GugPBS12,PolS10}.

\begin{remark} \label{rem:dae}
{
{\rm The  {pH} system class has been extended (e.g., in \cite{BeaMXZ18,MehM19,MehU23}),  to include {\em input-state-output {pH} descriptor ({pH}DAE) systems}, which take on a similar form (at least in the case of autonomous systems):
 \begin{equation} \label{{pH}daedef}
\begin{aligned}
   \bE\dot \bx&=\left(\bJ-\bR\right) \be(\bx)+(\bF-\bP)\bu(t), \\
\by(t)&=\ (\bF+\bP)^T \be(\bx)  + (\bS+\bN) \bu(t)
\end{aligned} \quad\mbox{ where }\nabla_{\!\bx}{\mathcal H}(\bx)=\bE^T \be(\bx),
\end{equation}
$\bJ$, $\bR$, $\bF$, $\bP$, $\bS$, and $\bN$ are as defined above, an auxiliary vector function $\be: \bRn \to \bRn$ has been introduced; and now $\bE=\bE^T\in \mathbb R^{n\times n}$ may be \emph{singular} (allowing then for the incorporation of algebraic constraints).  In this paper, we focus  on the case that $\bE$ is nonsingular and so without loss of generality it may be taken to be
the identity. }
}
\end{remark}
\pagebreak
{
\subsection{Transforming general LTI systems to pH form}
\label{sec:trafo}}
Consider now a general linear time-invariant (LTI) system:
\begin{equation} \label{GenSys}
\begin{array}{c}
\dot{\bx}\ =\ \bA\,\bx + \bB\, \bu,\\[1mm]
\by\ =\ \bC\,\bx + \bD\, \bu,
\end{array}
\end{equation}
with $\bA\in  \bRnn$, $\bB\in\bRnm$, $\bC\in\bRmn$, and $\bD \in\bRmm$.
Following the previous discussion leading to \eqref{DissipIneq}, the system (\ref{GenSys}) is \emph{passive} if there exists a continuously differentiable \emph{storage function}
${\mathcal H}:\bRn\rightarrow [0,\infty)$ such that (\ref{DissipIneq}) holds for all admissible inputs $\bu$ \cite{Wil72b}.

{ A natural pair of questions arise: \emph{When is (\ref{GenSys}) equivalent (in a sense made precise below) to a port-Hamiltonian (descriptor) system?  How may one construct an associated equivalence transformation,  implemented as numerically reliable procedure?}}

{ It is well-known \cite{BeaMV19,FauMPSW22} how one can use  \emph{equivalence transformations} to transform general LTI systems to port-Hamiltonian systems having the form}
\begin{equation} \label{FullpHSys}
\begin{array}{rcl}
\dot{\bfxi}&=&(\bJ-\bR)\bQ\, \bfxi\ +\ (\bF-\bP)\, \bfphi,\\[2mm]
\bfeta&=&(\bF+\bP)^\top\bQ\, \bfxi\  +\ (\bS+\bN)\, \bfphi,
\end{array}
\end{equation}
with  $\bJ=-\bJ^\top$, $\bR= \bR^\top\geq 0$, $\bQ= \bQ^\top>0$, $\bS=\bS^\top\ge 0$,
$\bN=-\bN^\top$, where $\bJ,\,\bR,\,\bQ\in {\mathbb R}^{n\times n}$,  $\bF,\bP\in {\mathbb R}^{n\times m}$,
$\bS,\bN\in {\mathbb R}^{m\times m}$, and $\bK$ as defined in (\ref{Kdef}) is positive semidefinite.

\begin{remark}\label{rem:QE}{\rm Note that by introducing $\bfzeta=\bQ\bfxi$ and $\bE=\bQ^{-1}$, we could alternatively discuss an equivalent descriptor formulation
\begin{equation} \label{FulldaepHSys}
\begin{array}{rcl}
\bE \dot{\bfzeta}&=&(\bJ-\bR)\, \bfzeta\ +\ (\bF-\bP)\, \bfphi,\\[2mm]
\bfeta&=&(\bF+\bP)^\top\, \bfzeta\  +\ (\bS+\bN)\, \bfphi.
\end{array}
\end{equation}}
\end{remark}

{ We will employ the following notion of \emph{system equivalence}, focussing on three invertible
}
transformations connecting (\ref{GenSys}) and (\ref{FullpHSys}), one on each of the input, the output, and the state space:
\[
\bu(t)= \widetilde{\bV} \bfphi(t), \  \bfeta(t) =\bV^\top \by(t),\ \mbox{and}\ \bx(t)=\bT^{-1} \bfxi(t) \
(\mbox{with }\widetilde{\bV}, \ \bV, \ \bT \ \mbox{invertible}).
\]
Within this context, the \emph{supply rate} associated with (\ref{GenSys}) is transformed as
\[
\by(t)^\top\bu(t)= \bfeta(t)^\top\bV^{-1}\widetilde{\bV}\bfphi(t).
\]
We wish to constrain the permissible transformations characterizing \emph{system equivalence}
so as to be \emph{power conserving};
that is, so that supply rates remain invariant, i.e. $ \by(t)^\top\bu(t)=\bfeta(t)^\top\bfphi(t)$. To guarantee this, we assume that $\widetilde{\bV}=\bV$ and we say that
(\ref{GenSys}) is \emph{equivalent} to a system of the form (\ref{FullpHSys}) if there exist invertible
matrices, $\bV$ and $\bT$, such that
\begin{equation} \label{IOtransformations}
\bu(t)=\bV\, \bfphi(t), \quad  \bfeta(t) =\bV^\top\, \by(t),\quad\mbox{and}\quad \bx(t)=\bT^{-1}\, \bfxi(t),
\end{equation}
and
\begin{eqnarray*}
&&\tilde \bA=\bJ-\bR= \bT \bA \bT^{-1},\quad \tilde \bB= \bF-\bP=\bT \bB \bV,\\
&&\tilde \bC=(\bF+\bP)^\top= \bV^\top\bC \bT^{-1},\quad \tilde \bD=\bS+\bN =\bV^\top\bD \bV,\quad \bQ=\bI.
\end{eqnarray*}
{ Although $\bV$ and $\bT$ need only be invertible to be candidate transformations, the freedom we have in choosing a final pH realization allows for $\bV$ to be assumed orthogonal and $\bT$ to be constructed from compositions of orthogonal and well-conditioned triangular transformations.}

\section{Passive systems and pH realizations}\label{sec:pa_and_ph}
{ In this section we first recall the construction of pH realizations for minimal systems and then extend this construction to the case of non-minimal systems.
\subsection{pH realization of minimal passive systems}
 Since pH systems are structurally passive \cite{CheGH22}, our starting point for the construction of transformations to pH form} is the following
characterization of passivity  introduced in \cite{Wil72b} for minimal linear time invariant systems.
The system (\ref{GenSys}) is \emph{minimal} if it is both controllable and observable.
The system (\ref{GenSys}) (and more specifically, the pair of matrices $(\bA,\bB)$ with $\bA \in \mathbb R^{n\times n}$, $\bB \in \mathbb R^{n\times m}$) is \emph{controllable} if $\rank \mat{cc} s\bI-\bA &\bB \rix=n$ for all $s\in \mathbb{C}$.
Similarly,  the system (\ref{GenSys}) (and the pair $(\bA,\bC)$ with $\bA \in \mathbb R^{n\times n}$, $\bC \in \mathbb R^{m\times n}$) is \emph{observable} if $\rank\mat{c}s\bI-\bA\\\bC\rix=n$ for all $s\in \mathbb{C}$.

\begin{theorem}[\cite{Wil72b}]\label{Willems1972}
Assume that the LTI system (\ref{GenSys}) is minimal.
The {Kalman-Yakubovich-Popov (KYP) linear} matrix inequality
\begin{equation} \label{LMIcond}
\left[\begin{array}{cr} \bA^\top\bQ +\bQ \bA & \bQ\bB- \bC^\top \\[2mm]
\bB^\top\bQ -\bC  & -(\bD+\bD^\top)
\end{array}\right]\leq 0
\end{equation}
has a solution $\bQ= \bQ^\top> 0$ if and only if (\ref{GenSys}) is a \emph{passive} system, in which case:
\begin{itemize}
\item [i)] ${\mathcal H}(\bx) = \frac12 \bx^\top \bQ \bx$ defines a storage function for (\ref{GenSys}) associated
with the supply rate $\by^\top\bu$, satisfying (\ref{DissipIneq}).
\item [ii)] There exist  maximum and minimum symmetric solutions to (\ref{LMIcond}): $\bQ_{+}\geq \bQ_{-} > 0$ such that
for all symmetric solutions $\bQ$ to (\ref{LMIcond}), $\bQ_{-}\leq\bQ\leq  \bQ_{+}$.
\end{itemize}
\end{theorem}
{ Using the transformations in  \eqref{IOtransformations}, this result has an immediate consequence for \emph{pH realizations}, i.e. equivalent representations in pH form, see also \cite{CheGH22}.}
 \begin{corollary}\label{minimalPassive}
Assume that the LTI system (\ref{GenSys}) is minimal.  Then
(\ref{GenSys}) has a pH realization ({ i.e.,
(\ref{GenSys}) is equivalent to (\ref{FullpHSys})}) if and only if it is passive.  Moreover, if (\ref{GenSys}) is passive then \emph{every} system equivalent to (\ref{GenSys})  (as generated by transformations in (\ref{IOtransformations})) is directly expressible as a pH system of the form \eqref{FullpHSys}.
\end{corollary}
An explicit numerical construction of a pH realization of a passive system can be performed as follows: {
If (\ref{GenSys}) is passive then (\ref{LMIcond}) has a positive definite solution
$\widehat{\bQ}= \widehat{\bQ}^\top=\bT^\top\bT$, (e.g., written in terms of a Cholesky factorization)}.  Then we can define directly
\begin{equation}\label{pHelem}
\begin{array}{cll}
\bQ=\bI, &  \bJ=\frac12(\bT \bA\bT^{-1}-(\bT\bA\bT^{-1})^\top), & \bR=-\frac12(\bT \bA\bT^{-1}+(\bT\bA\bT^{-1})^\top)\\[1mm]
& \bF=\frac12\left(\bT \bB+(\bC\bT^{-1})^\top\right), & \bP=\frac12\left(-\bT \bB+(\bC\bT^{-1})^\top\right),\\[1mm]
& \bS=\frac12( \bD+\bD^\top) & \bN=\frac12( \bD-\bD^\top).
\end{array}
\end{equation}
Since  (\ref{LMIcond}) can be written in terms of these defined quantities as
\[
-2\ \left[\begin{array}{cc}
\bT^\top & 0 \\
0 & \bI
\end{array}\right]
\ \left[\begin{array}{lc}
\bR & \bP \\
\bP^\top & \bS
\end{array}\right]\ \left[\begin{array}{cc}
\bT & 0 \\
0 & \bI
\end{array}\right] \ \leq 0,
\]
this gives  (\ref{Kdef}) and
$\bJ,\,\bR,\,\bQ,\,\bF,\,\bP,\,\bS,\,\mbox{and } \bN$
as defined in (\ref{pHelem}) will indeed determine a pH system.

\begin{remark}\label{rem:des}{\rm
Note that instead of the transformation in \eqref{pHelem} we may use the descriptor formulation \eqref{FulldaepHSys}
with
\begin{equation}\label{pHdaeelem}
\begin{array}{cll}
\bE=\widehat{\bQ}, &  \bJ=\frac12( \widehat{\bQ}\bA-(\widehat{\bQ}\bA)^\top), & \bR=-\frac12(\widehat{\bQ}\bA+(\widehat{\bQ}\bA)^\top),\\[1mm]
& \bF=\frac12\left(\widehat{\bQ}\bB+\bC^\top\right), & \bP=\frac12\left(-\widehat{\bQ}\bB+\bC^\top\right),\\[1mm]
& \bS=\frac12( \bD+\bD^\top), & \bN=\frac12( \bD-\bD^\top)
\end{array}
\end{equation}
that avoids the factorization of $\widehat{\bQ}$ and the similarity transformation with $\bT$.
}
\end{remark}
The construction described above for  minimal passive systems is well known,
{ see e.g. \cite{BeaMV19}. In the next subsection we discuss how to do this for non-minimal passive systems.
\subsection{Construction for non-minimal systems}\label{sec:nonminimal}}
In the previous subsection we  have discussed the existence of a positive
definite solution $\bQ$ of \eqref{LMIcond} under the assumption of minimality of the system. However, such solutions may exist even if the system is not minimal.  A detailed analysis of the intricate relationship between passive systems, port-Hamiltonian descriptor systems and the solvability of the KYP inequality has recently been presented in \cite{CheGH22} and an explicit construction is presented in \cite{WeiWS94}. This analysis is particularly important when the system matrices arise from an interpolatory realization or a model reduction process where the resulting systems may be non-minimal or otherwise very close approximations of non-minimal systems. In these circumstances, the computation of a pH representation may be very sensitive to small perturbations arising from measurement or round-off errors.

Consider the following example  from \cite{BeaMV19}.
\begin{example}\label{rem:illcondex}{\rm
The system
$\dot{x}=-x,y = u$  is both stable and passive but not minimal.
In this case, the inequality  (\ref{LMIcond}) is satisfied with any (scalar)
{
$Q>0$, and $J=0$, $R=1/Q$, $B=C=0$, and $D=1$. The Hamiltonian may be defined as
$\mathcal H(x)=\frac{Q}{2} x(t)^2$, and the dissipation inequality evidently holds since, for $t_1\geq t_0$,
\begin{align*}
\mathcal{H}(x(t_1))-\mathcal{H}(x(t_0)) & = \frac{Q}{2} (x(t_0) e^{-(t_1-t_0)})^2 -\frac{Q}{2} x(t_0)^2 \\
 = \frac{Q}{2} x(t_0)^2& (e^{-2(t_1-t_0)}-1) \leq 0 \leq \int_{t_0}^{t_1} y(t) \, u(t) \, dt = \int_{t_0}^{t_1} u(t)^2 \, dt.
\end{align*} }
}
\end{example}
{In Section~\ref{sec:Ricin}, we analyze how the
conditions for the existence of solutions can be verified numerically  and how one may calculate solutions to \eqref{LMIcond}.}

{Note that when a system of the form \eqref{GenSys} is generated by an interpolatory realization or other model reduction strategies, it may only be a close approximation to a passive system even when the original system is passive. In this case the inequality~\eqref{LMIcond} might not be solvable, but a solution may exist for an adjacent system obtained from a small perturbation of the coefficients $\bA,\bB,\bC,\bD$. How to obtain the ``best" such perturbation is an important research topic, see \cite{AlaBKMM11,BruS13,Gri04} and references therein. We develop such perturbations when needed in the course of our construction, and so (referencing the pH representation), we are able to replace a nearly passive system with a nearby passive one.}

\begin{remark}\label{rem:freedom}
 {\rm Since the matrix inequality \eqref{LMIcond} typically has an infinite number of solutions, an important question  is how best to use the freedom in the choice of the solution of \eqref{LMIcond}. One natural goal might be to minimize the distance to either instability or non-passivity or to maximize other robustness measures by pursuing the so-called analytic center of the solution set, see \cite{BanMNV20,MehV20}.  Characterizing the solution set of \eqref{LMIcond} and its relationship to different robustness measures remains currently a (mostly) open problem, but see \cite{MehX24} for recent partial results.}
\end{remark}

The matrix inequality (\ref{LMIcond}) implies
the \emph{Lyapunov inequality},
\begin{equation}\label{lyain}
\bA^\top\bQ +\bQ \bA\leq 0,
\end{equation}
and via Lyapunov's theorem \cite{LanT85}, this guarantees in turn that the unforced system $\dot \bx=\bA \bx$ is \emph{stable}, i.e., $\bA$ has all eigenvalues in the closed left half plane and those on the imaginary axis are semisimple;
if the inequality (\ref{lyain}) is strict, then the system is \emph{asymptotically stable}, i.e., all eigenvalues of $\bA$ are in the open left half plane.
Passivity is encoded in the solvability of the full matrix inequality~(\ref{LMIcond}); \emph{strict passivity} occurs when the inequality is strict, i.e., if the dissipation inequality (\ref{DissipIneq}) is strict, see \cite{KotA10} for a detailed analysis.

\begin{remark}\label{rem:extreme_case}{\rm
In order to characterize the boundary of the solution sets  of the LMIs~(\ref{LMIcond}) and (\ref{lyain}) one needs to study the case when either the inequalities in (\ref{LMIcond}) and (\ref{lyain}) are not strict or when the resulting solutions are only semidefinite (or both).
Extreme points of the solution set of \eqref{LMIcond} have recently been characterized in \cite{MehX24}.
}
\end{remark}

{
This paper is organized as follows. In Section~\ref{sec:lyaric} we discuss the solution of the Lyapunov and Riccati matrix inequalities in the  general situation of non-minimal systems  and for the case that the symmetric part of $\bD$ is not invertible in Section~\ref{sec:suff}. We recall previous results and  present numerical procedures to  construct
explicit transformations mapping a general linear time-invariant system  to a port-Hamiltonian form (\ref{FullpHSys}) in  Section~\ref{sec:numerical}.}
\section{Lyapunov and Riccati inequalities}\label{sec:lyaric}

The solutions of Lyapunov and Riccati inequalities as they arise in
the characterizations (\ref{LMIcond}) and  (\ref{lyain}) are typically addressed through semidefinite programming, see \cite{BoyEFB94}. { For (\ref{LMIcond}), an explicit characterization
for all possible linear system realizations can be found in \cite{AlpL11}. 
In \cite{WeiWS94} a constructive method was proposed to deal with the case when $\bD+\bD^\top$ is singular. In contrast to \cite{AlpL11}, we
assume that a general linear state-space system is given. We derive a constructive step-by-step procedure to check  whether the linear system is passive,
and if so, characterize all possible positive definite solutions of (\ref{LMIcond}) for transforming (\ref{GenSys}) to a pH system (\ref{FullpHSys}).
}

\subsection{Solution of Lyapunov inequalities.}
The stability of $ \bA$ is a necessary condition for (\ref{lyain}) and (\ref{LMIcond})
which require that $\bT \bA\bT^{-1}+\bT^{-\top} \bA^\top\bT^\top\leq 0$, or equivalently,  that the Lyapunov inequality (\ref{lyain})
has a positive definite solution $\bQ=\bT^\top\bT$.
It is well known that the equality case in (\ref{lyain}) always has a positive definite solution if $ \bA$ is stable,
see \cite{LanT85}. In the following we recall, see e.g. \cite{BoyEFB94}, a characterization of  the complete set of solutions of the inequality case.

If $ \bA$ is stable, but not asymptotically stable, then since all eigenvalues on the imaginary axis must be semi-simple,   the real Jordan form of $ \bA$, see e.g. \cite{HorJ85}, guarantees the existance of a nonsingular matrix $\bM\in \bRnn$ such that
\eq{decA}
\bM \bA\bM^{-1}=\diag\left(\bA_1,\alpha_2\bJ_2,\ldots,\alpha_r\bJ_r\right),
\en
where $\bA_1\in \mathbb R^{n_1\times n_1}$ is asymptotically stable, $\alpha_2,\ldots,\alpha_{r}\ge 0$ are real
and distinct, and $\bJ_j=\mat{cc}0&\bI_{n_j}\\-\bI_{n_j}&0\rix$, $j=2,\ldots,r$. In order to characterize the solution set of  (\ref{lyain}), we make the ansatz
\eq{forX}
\bQ=\bM^\top\diag\left(\bQ_1,\hat \bQ_2,\ldots,\hat \bQ_{r}\right)\bM,
\en
and  consider the determination of the block $\bQ_1$ separately from determination of the other blocks.
Let  $\mathcal W(n_1)$ be the set
of symmetric positive semidefinite matrices $\Theta_1\in \mathbb R^{n_1\times n_1}$ with the property that $\Theta_1\bx\ne 0$ for any eigenvector  $\bx$ of $\bA_1$. Then for any $\Theta_1\in \mathcal W(n_1)$ we define $\bQ_1$ to be the unique symmetric positive definite  solution of the Lyapunov equation $\bA_1^\top\bQ_1+\bQ_1\bA_1=-\Theta_1$, see \cite{LanT85}.
The other matrices $\hat \bQ_j$, $j=2,\ldots,r$ are chosen of the form
\[
\hat \bQ_j=\mat{cc}\bY_j&\bZ_j\\-\bZ_j&\bY_j\rix>0,
\]
with $\bZ_j=-\bZ_j^\top$, when $\alpha_{j}>0$ or an arbitrary $\hat \bQ_j>0$ when $\alpha_{j}=0$.

We have the following characterization of the solution set of (\ref{lyain}), {see \cite[Section 2.5.2]{BoyEFB94}}.
\begin{lemma}\label{LyapStab}
Let $ \bA\in \bRnn$. Then the Lyapunov inequality (\ref{lyain})
has a symmetric
positive definite solution $\bQ\in \bRnn$ if and only if $\bA$ is stable.

If $\bA$ is asymptotically stable, then the solution set is given by the set of all symmetric positive definite solutions of the Lyapunov equation $ \bA^\top\bQ+\bQ\bA=-\Theta$, where $\Theta$ is any symmetric positive semidefinite matrix  with the property that $\Theta\bx\ne 0$ for any eigenvector  $\bx$ of $ \bA$, i.e.,
$( \bA,\Theta)$ is observable.

If $\bA$ is stable, but not asymptotically stable, then with the transformation (\ref{decA}), any solution of (\ref{lyain}) must have the form (\ref{forX}), solving the Lyapunov equation
\eq{forW}
\bA^\top\bQ+\bQ \bA=-\bM^\top\diag(\Theta_1,0,\ldots,0)\bM
\en
with $\Theta_1\in \mathcal W(n_1)$.
\end{lemma}

\begin{remark}{\rm 
Lemma~\ref{LyapStab}
relies on the computation of the real Jordan form of $ \bA$,
which generally is problematic in finite precision arithmetic.
For the numerical computation of the solution it is better to use the real Schur form,
see \cite{GolV96}.
}
\end{remark}
{
\subsection{Solution of the KYP inequality
in the case $\bD+\bD^\top>0$.}
\label{sec:ricineq}}
Using Corollary~\ref{minimalPassive} we can characterize (at least in the minimal case) the existence of a transformation to pH form
via  the existence of a symmetric positive definite matrix $\bQ$ solving the {KYP linear matrix inequality}
{\color{black}
\eq{LMIred}
\mat{cc}\bA^\top\bQ+\bQ\bA&\bQ\bB-\bC^\top\\
\bB^\top\bQ-\bC&-(\bD+ \bD^\top)\rix\leq 0.
\en
}
{
In Section~\ref{sec:suff} we  will develop a numerical method to
check whether the KYP linear matrix inequality (\ref{LMIred}) has a positive definite solution, and we characterize all possible positive definite solutions when they do exist. In this section we first assume that $\hat \bS:=\bD+\bD^\top>0$ and recall results for this case. The singular case will be treated in Section~\ref{singularS}.}
%

We extend the well-known results about the solvability of \eqref{LMIred} to the general case that the system may be either non-controllable or non-observable.  For this and for the numerical methods, we will have to identify the controllable and observable subsystems in a numerically viable way.

If $\hat \bS>0$, by using Schur complements, we have that  (\ref{LMIred}) is equivalent to the \emph{Riccati inequality}
\begin{equation}\label{RIE}
( \bA- \bB \hat \bS^{-1} \bC)^\top \bQ
+ \bQ( \bA- \bB  \hat \bS^{-1} \bC)
+\bQ \bB \hat \bS^{-1} \bB^\top\bQ+  \bC^\top  \hat \bS^{-1}  \bC\le 0.
\end{equation}
To study the solvability of \eqref{RIE}, we first investigate the influence of the purely imaginary eigenvalues of $\bA$ on the solvability of (\ref{LMIred}) and \eqref{RIE}.
Clearly, a necessary condition for (\ref{LMIred}) to be solvable is that
$\bQ$ satisfies $\bA^\top\bQ+\bQ\bA\le 0$. Following Lemma~\ref{LyapStab},
if $\bA$ has the form (\ref{decA}) then $\bQ$ must have the form
(\ref{forX}), and $ \bA^\top\bQ+\bQ\bA$ has the form (\ref{forW}). Written in
compact form, we get
\[
\bM\bA\bM^{-1} =\mat{cc}\bA_1&0\\0&\bA_2\rix,\quad
\bM^{-\top}\bQ\bM^{-1}=\mat{cc}\bQ_1&0\\0&\bQ_2\rix,
\]
where
\[
 \bA_2=\diag(\alpha_2\bJ_2,\ldots,\alpha_r\bJ_r)=-\bA_2^\top,\quad
\bQ_2 =\diag(\hat \bQ_2,\ldots,\hat \bQ_r),
\]
satisfies $\bA_2^\top\bQ_2+\bQ_2\bA_2=0$, and
\[
\bM^{-\top}(\bA^\top\bQ+\bQ\bA)\bM^{-1}
=\mat{cc}\bA_1^\top\bQ_1+\bQ_1\bA_1&0\\0&0\rix.
\]
Setting
\[
\bM\bB =\mat{c}\bB_1\\\bB_2\rix,\quad
\bC\bM^{-1}=\mat{cc}\bC_1&\bC_2\rix.
\]
and premultiplying $\bM^{-\top}$ and post-multiplying $\bM^{-1}$ to the
first block row and column of (\ref{LMIred}), respectively, one has
that
{\color{black}
\[
\mat{ccc} \bA_1^\top\bQ_1+\bQ_1\bA_1&0&\bQ_1\bB_1-\bC_1^\top\\
0&0&\bQ_2\bB_2-\bC_2^\top\\
\bB_1^\top\bQ_1-\bC_1&\bB_2^\top\bQ_2-\bC_2& -\hat \bS\rix\le 0.
\]
}
Therefore, to have a positive definite solution of \eqref{LMIcond}, $\bQ_2$ must be positive definite satisfying
\eq{nec1}
\bB_2^\top\bQ_2=\bC_2,\quad \bA_2^\top\bQ_2+\bQ_2\bA_2=0,
\en
and $\bQ_1$ must be a positive definite solution of the smaller size linear  matrix inequality
{\eq{LMIredred}
\mat{cc}\bA_1^\top\bQ_1+\bQ_1\bA_1&\bQ_1\bB_1-\bC_1^\top\\
\bB_1^\top\bQ_1-\bC_1&-\hat \bS\rix\le 0
\en
}
or equivalently $\bQ_1$ has to satisfy the \emph{ Riccati inequality}
\begin{eqnarray}\nonumber
\Psi(\bQ_1)&:=&
( \bA_1- \bB_1 \hat \bS^{-1} \bC_1)^\top\bQ_1
+ \bQ_1( \bA_1- \bB_1  \hat \bS^{-1} \bC_1)\\
&+&\bQ_1 \bB_1 \hat \bS^{-1} \bB_1^\top\bQ_1
+  \bC_1^\top  \hat \bS^{-1}  \bC_1\le 0.\label{RIEred}
\end{eqnarray}
To study the solvability of (\ref{RIEred}), we recall that by construction, $\bA_1$ is asymptotically stable.
We claim that { $\bA_1-\bB_1\hat \bS^{-1}\bC_1$} is necessarily asymptotically
stable as well. To show this, suppose that the inequality (\ref{RIEred}) has a solution $\bQ_1>0$. Since
$\bQ_1 \bB_1 \hat \bS^{-1} \bB_1^\top\bQ_1+  \bC_1^\top
 \hat \bS^{-1}  \bC_1\ge 0$,
it follows from Lemma~\ref{LyapStab}  that { $ \bA_1- \bB_1 \hat \bS^{-1} \bC_1$} is stable, and
 there exists an invertible matrix $\bM_1$ such that
$\bM_1(\bA_1-\bB_1\hat \bS^{-1}\bC_1)\bM_1^{-1}=\diag(\widetilde \bA_1,\alpha_2\bJ_2,\ldots,\alpha_r\bJ_r)$ is in real Jordan form as in
(\ref{decA}), where $\wt\bA_1$ is asymptotically stable, $\bM_1^{-\top}\bQ_1\bM_1^{-1}$ has the form (\ref{forX}) and following
(\ref{forW}), we have
\[
\bM_1^{-\top}(\bQ_1 \bB_1 \hat \bS^{-1} \bB_1^\top\bQ_1
+  \bC_1^\top  \hat \bS^{-1}  \bC_1) \bM_1^{-1}
=\diag(\widetilde\Theta,0,\ldots,0)\ge 0.
\]
Then, due to the positive definiteness of $ \hat \bS$ it follows that
$ \bC_1\bM_1^{-1}=\mat{cccc}\bC_{11}&0&\ldots&0\rix$,
and by making use of the block diagonal structure of
$\bM_1^{-\top}\bQ_1\bM_1^{-1}$, 
we also have
$\bM_1\bB_1=\mat{cccc}\bB_{11}^\top&0&\ldots&0\rix^\top$.
Thus it follows that
\[
\bM_1 \bA_1\bM_1^{-1}
=\diag({\wt\bA_{1}}+\bB_{11}\hat \bS^{-1}\bC_{11},\alpha_2\bJ_2,\ldots,\alpha_r\bJ_r).
\]
Since $\bA_1$ is asymptotically stable, all $\alpha_j\bJ_j$ must be void, which
implies that $\bA_1-\bB_1\hat \bS^{-1}\bC_1$ must be
asymptotically stable as well.

In order to characterize the solution of the  Riccati inequality
(\ref{RIEred}), we first have to identify what happens if the system is not minimal. To numerically check minimality, we can use the orthogonal version of the
Kalman decomposition, \cite{Kai80,Van79}, see also \cite{MehX00}.
\begin{lemma}\label{lem:stcf} Consider a  general system of the form (\ref{GenSys}).
Then there exists a real orthogonal matrix $ \bU$ such that
\begin{eqnarray}
\nonumber
 \bU^\top\bA \bU
&=&
\mat{cc|c}\widehat \bA_{11}&0&\widehat \bA_{13}\\\widehat \bA_{21}&\widehat \bA_{22}&\widehat \bA_{23}\\
\hline
0&0&\widehat \bA_{33}\rix
=:\mat{cc}\wt \bA_{11}&\wt\bA_{12}\\0&\wt\bA_{22}\rix\\
&&\label{stcf}\\
\nonumber
 \bU^\top\bB&=&\mat{c}\widehat \bB_{1}\\ \widehat \bB_{2}\\\hline 0\rix=:
\mat{c}\wt\bB_1\\0\rix,\quad
\bC \bU=\mat{cc|cc}\widehat \bC_{1}&0&\widehat \bC_{3}\rix
=:\mat{cc}\wt\bC_1&\wt\bC_2\rix,
\end{eqnarray}
where the pairs
$(\wt\bA_{11},\wt\bB_1)$ and $(\widehat \bA_{11},\widehat \bB_{1})$ are controllable and  the pair $(\widehat \bA_{11},\widehat \bC_{1})$ is observable.
\end{lemma}

The next lemma considers the Riccati inequality  (\ref{RIEred}), where now the coefficients are transformed to the form in \eqref{stcf}.
\begin{lemma} For the Riccati inequality (\ref{RIEred}), suppose  that both $ \bA_1$
and $\bA_1-\bB_1\hat \bS^{-1}\bC_1$ are asymptotically stable, and that
 $\hat \bS>0$.
Then there exists a transformation to the condensed form (\ref{stcf}) of $\bA_1,\bB_1,\bC_1$, such that
\eqn
\nonumber
  \bU^\top( \bA_1- \bB_1 \hat \bS^{-1}  \bC_1) \bU
&=&\mat{cc}\wt\bA_{11}-\wt\bB_1\hat \bS^{-1}\wt\bC_1
&\wt\bA_{12}-\wt\bB_1\hat \bS^{-1}\wt\bC_2\\
0&\wt\bA_{22}\rix
\\
\label{defa1}
&=&\mat{cc|c}\bA_{11}-\bB_{1}\hat \bS^{-1}\bC_{1}&0&\bA_{13}-\bB_{1}\hat \bS^{-1}\bC_{3}
\\\bA_{21}-\bB_{2}\hat \bS^{-1}\bC_{1}&\bA_{22}&\bA_{23}-\bB_{2}\hat \bS^{-1}\bC_{3}\\\hline
0&0&\bA_{33}\rix,
\enn
where $\bA_{11}-\bB_{1}\hat \bS^{-1}\bC_{1}$, $\bA_{22}$, and $\bA_{33}$ are all
asymptotically stable. In addition,  $\bA_{11}$  also is asymptotically stable, $(\wt\bA_{11},\wt\bB_{1})$
and $(\bA_{11},\bB_1)$ are controllable, and $(\bA_{11},\bC_{1})$ is  observable.
\end{lemma}
\proof
The proof is straightforward and omitted.
\eproof

Using the previous lemma we may check controllability and observability, {though it is well-known \cite{PaiV81} that these properties} can also be read off from Lagrangian invariant subspaces (if they exist) of certain Hamiltonian matrices associated with the equality case in \eqref{RIE}, respectively \eqref{RIEred}, see e.g. \cite{LanR95}.
\begin{lemma}~\label{lem:CS}
Suppose that   $\hat \bS>0$ and that  the {\em Hamiltonian matrix}
\begin{equation}\label{hampas}
\bH
:=\mat{cc} \bA-\bB\hat \bS^{-1}\bC&\bB\hat \bS^{-1}\bB^\top \\
-\bC^\top\hat \bS^{-1} \bC &  -(\bA-\bB\hat \bS^{-1}\bC)^\top
\rix
\end{equation}
has a Lagrangian invariant subspace, i.e., there exist square matrices
$\bW_1$, $\bW_2$, and $\bZ$ such that
\begin{equation}\label{haminv}
    \bW_1^T\bW_2=\bW_2^T\bW_1,\
    \begin{bmatrix}
        \bW_1 \\ \bW_2
    \end{bmatrix}\mbox{ has full column rank, and }    \bH \begin{bmatrix}
        \bW_1 \\ \bW_2
        \end{bmatrix}
            =
        \begin{bmatrix}
        \bW_1 \\ \bW_2
        \end{bmatrix} \bZ.
\end{equation}
Then, if the pair $(\bA,\bB)$ is controllable, $\bW_1$ is invertible;  if $(\bA,\bC)$ is observable then  $\bW_2$ is invertible.
\end{lemma}

Note that { in Lemma~\ref{lem:CS} the matrix} $\bW_1$ may be still invertible even if $(\bA,\bB)$ is
not controllable, for instance, when $\bB=0$ and $\bC=0$. The same
applies to $\bW_2$.

The following lemma is also well-known, see e.g. \cite{LanR95,Meh91}.
\begin{lemma}~\label{lem:COSol}
Suppose that $(\bA,\bB)$ is controllable, $(\bA,\bC)$
is observable, $\bA$ is asymptotically
stable, and $\hat \bS>0$. Then the Riccati equation
\eq{ARE}
(\bA-\bB\hat \bS^{-1}\bC)^\top\bQ+\bQ
(\bA-\hat \bS^{-1}\bC)+\bQ\bB\hat \bS^{-1}\bB^\top\bQ
+\bC^\top\hat \bS^{-1}\bC=0
\en
has a solution $\bQ>0$ if and only if the Hamiltonian matrix in \eqref{hampas}
has a Lagrangian invariant subspace satisfying \eqref{haminv}. If such a Lagrangian
invariant subspace exists, then $\bQ:=\bW_2\bW_1^{-1}>0$
solves (\ref{ARE}).
\end{lemma}

\begin{remark}\label{rem:notobserv}{\rm
By Lemma~\ref{lem:CS} we have seen that if the pair $(\bA,\bC)$ is not observable then the matrix $\bW_2$ in \eqref{haminv}
may or may not be invertible. Hence, if $\bW_1$ is invertible and $\bW_2$ is
not, then a symmetric solution to \eqref{ARE} still exits but is only positive semidefinite. This gives a characterization of the boundary of the solution set of the matrix inequality~\eqref{LMIred}, see also \cite{CheGH22} for a detailed discussion of  the existence of solution to the matrix inequality \eqref{LMIred} if the controllability and observability conditions are violated.
}
\end{remark}


Lemma~\ref{lem:COSol} shows that under the conditions of controllability and observability the Riccati equation (\ref{ARE}) has a
solution $\bQ>0$ whenever the Hamiltonian matrix $\bH$ has
a Lagrangian invariant subspace. The existence of such an
invariant subspace depends only on the  purely
imaginary eigenvalues of $\bH$, e.g., \cite{FreMX02}. When
such an invariant subspace exists, then there are many such
invariant subspaces. Note that the eigenvalues of $\bH$ are
symmetric with respect to the imaginary axis in the complex plane. If
(\ref{haminv}) holds, then the union of the
eigenvalues of $\bZ$ and $-\bZ^\top$ form the spectrum of $\bH$.
One particular choice is that the spectrum of $\bZ$ is in the closed left
half complex plane, another choice is that it is in the closed right half
complex plane. The two corresponding solutions of
the Riccati equation (\ref{ARE})
are the {\em minimal} solution $\bQ_-$ and the {\em maximal} solution $\bQ_+$ and all other solutions of the Riccati equation lie (in the Loewner ordering of symmetric matrices) between these extremal solutions.

\begin{example}{\rm
Consider the example $\bB=\hat \bS=1$, $\bC=-1$, $\bA=-1-\alpha$,
where $\alpha>0$. So $\bA$ and $\bA-\bB\hat \bS^{-1}\bC=-\alpha$
are both asymptotically stable.
The Riccati equation (\ref{ARE}) is
$q^2-2\alpha q+1=0$,
which does not have a real positive semidefinite solution when $\alpha\in (0,\,1)$.
If $\alpha=1$, then it has a unique solution $q=1$ associated with $\bZ=0$.
If $\alpha>1$, then it has two solutions $\alpha\pm\sqrt{\alpha^2-1}>0$
with $\bZ=\pm\sqrt{\alpha^2-1}$.
}
\end{example}

Suppose that $\bA_1,\bB_1,\bC_1$ from (\ref{RIEred})
have been transformed via an orthogonal matrix $\bU$ into the form (\ref{stcf}).
Partition
\[
\bU^\top\bQ_1\bU
=\mat{cc}\wt \bQ_{11}&\wt\bQ_{12}\\\wt\bQ_{12}^\top&\wt\bQ_{22}\rix.
\]
The  Riccati inequality (\ref{RIEred}) then is equivalent to
\eq{psiLMI}
\bU^\top \Psi(\bQ_1)\bU=\mat{cc}
\wt\Psi_{11}(\bQ_1)&\wt\Psi_{12}(\bQ_1)\\
\wt\Psi_{12}(\bQ_1)^\top&\wt\Psi_{22}(\bQ_1)\rix
\le 0,
\en
where (suppressing arguments for compactness)
\begin{eqnarray*}
\wt\Psi_{11} &=&(\wt \bA_{11}-\wt\bB_1\hat \bS^{-1}\wt\bC_1)^\top\wt\bQ_{11}+\wt\bQ_{11}(\wt\bA_{11}
-\wt\bB_1\hat \bS^{-1}\wt\bC_1)+
\wt\bQ_{11}\wt\bB_1\hat \bS^{-1}\wt\bB_1^\top\wt\bQ_{11}
+\wt\bC_1^\top\hat \bS^{-1}\wt\bC_1,\\
\wt\Psi_{12}&=&(\wt\bA_{11}-\wt\bB_1\hat \bS^{-1}(\wt\bC_1-
\wt\bB_1^\top\wt\bQ_{11}))^\top\wt\bQ_{12}+\wt\bQ_{12}\wt\bA_{22}
+\wt\bQ_{11}(\wt\bA_{12}-\wt\bB_1\hat \bS^{-1}\wt\bC_2)+
\wt\bC_1^\top\hat \bS^{-1}\wt\bC_2,\\
\wt\Psi_{22}&=&\wt\bA_{22}^\top\wt\bQ_{22}
+\wt\bQ_{22}\wt\bA_{22}
+(\wt\bA_{12}-\wt\bB_1\hat \bS^{-1}\wt\bC_2)^\top\wt\bQ_{12}+
\wt\bQ_{12}^\top(\wt\bA_{12}-\wt\bB_1\hat \bS^{-1}\wt\bC_2)\\
&&\qquad\qquad+
\wt\bQ_{12}^\top\wt\bB_1\bS^{-1}\wt\bB_1^\top\wt\bQ_{12}
+\wt\bC_2^\top\hat \bS^{-1}\wt\bC_2.
\end{eqnarray*}
For (\ref{RIEred}) to have a positive definite solution $\bQ_1>0$, it is necessary that
\[
\wt\Psi_{11}(\bQ_1)=\wt\Psi_{11}(\wt \bQ_{11})\le 0
\]
has a positive definite solution $\wt \bQ_{11}$
or equivalently, that the
dual Riccati inequality
\[
\Phi(\wt \bY):=
(\wt \bA_{11}-\wt\bB_1\hat \bS^{-1}\wt\bC_1)\wt \bY+
\wt\bY(\wt\bA_{11}-\wt\bB_1\hat \bS^{-1}\wt\bC_1)^\top+\wt \bY\wt\bC_1^\top\hat \bS^{-1}\wt\bC_1\wt\bY+\wt\bB_1\hat \bS^{-1}\wt\bB_1^\top\le 0
\]
has a solution $\wt\bY=\wt\bQ_{11}^{-1}>0$. Using the partitioning in (\ref{stcf})
for  $\wt\bY=\mat{cc}\bY_{11}&\bY_{12}\\\bY_{12}^\top&\bY_{22}\rix$, one can write this as
\[
\Phi(\wt\bY)=
\mat{cc}\Phi_{11}(\wt\bY)&\Phi_{12}(\wt\bY)\\\Phi_{12}(\wt\bY)^\top&\Phi_{22}
(\wt\bY)\rix\le 0,
\]
where
\begin{eqnarray*}
\Phi_{11}(\wt\bY)&=&\Phi_{11}(\bY_{11})=(\bA_{11}-\bB_1\hat \bS^{-1}\bC_1)\bY_{11}+\bY_{11}
(\bA_{11}-\bB_1\hat \bS^{-1}\bC_1)^\top\\ &+&\bY_{11}\bC_1^\top\hat \bS^{-1}\bC_1\bY_{11}
+\bB_1\hat \bS^{-1}\bB_1^\top.
\end{eqnarray*}
It is necessary that $\Phi_{11}(\bY_{11})\le 0$ has a solution $\bY_{11}>0$,
or equivalently that the dual inequality
\eq{RIE11}
(\bA_{11}-\bB_1\hat \bS^{-1}\bC_1)^\top\bQ_{11}+\bQ_{11}
(\bA_{11}-\bB_1\hat \bS^{-1}\bC_1)+\bQ_{11}\bB_1\hat \bS^{-1}\bB_1^\top\bQ_{11}
+\bC_1^\top\hat \bS^{-1}\bC_1\le 0
\en
has a solution $\bQ_{11}>0$. This is equivalent to the fact that the equality  case in \eqref{RIE11} has a
positive definite solution, see \cite{BoyEFB94,LanR95}.
Since the corresponding Hamiltonian matrix is
\eq{hatH}
\bH_{11}=\mat{cc}\bA_{11}-\bB_1\hat \bS^{-1}\bC_1&\bB_1\hat \bS^{-1}\bB_1^\top\\
-\bC_1^\top\hat \bS^{-1}\bC_1&-(\bA_{11}-\bB_1\hat \bS^{-1}\bC_1)^\top\rix,
\en
and since $(\bA_{11},\bB_1,\bC_1)$ is minimal, by Lemma~\ref{lem:COSol},
 \eqref{RIE11} has a solution $\bQ_{11}>0$ if and only if $\bH_{11}$ has a Lagrangian invariant subspace.
This and the condition that { $\bA_{11}-\bB_1\hat \bS^{-1}\bC_1$ }is asymptotically stable are necessary
conditions for the solvability of (\ref{LMIredred}) or (\ref{RIEred}).

We now show through an explicit construction that these two conditions are also sufficient
for the existence of a positive definite solution of (\ref{LMIredred}) or (\ref{RIEred}).
Together with \eqref{nec1}, they constitute necessary and sufficient conditions for
the solvability of (\ref{LMIred}) or (\ref{RIE}).

The Hamiltonian matrix corresponding to $\wt\Psi_{11}(\wt \bQ_{11})\le 0$ is
\begin{eqnarray*}
\wt \bH_{11}&=&\mat{cc}\wt\bA_{11}-\wt\bB_1\hat \bS^{-1}\wt\bC_1&\wt\bB_1\hat \bS^{-1}\wt\bB_1^\top\\
-\wt\bC_1^\top\hat \bS^{-1}\wt\bC_1&-(\bA_{11}-\wt\bB_1\hat \bS^{-1}\wt\bC_1)^\top\rix\\
&=&\mat{cc|cc}\bA_{11}-\bB_1\hat \bS^{-1}\bC_1&0&\bB_1\hat \bS^{-1}\bB_1^\top&
\bB_1\hat \bS^{-1}\bB_2^\top\\
\bA_{21}-\bB_2\hat \bS^{-1}\bC_1&\bA_{22}&\bB_2\hat \bS^{-1}\bB_1^\top&\bB_2\hat \bS^{-1}\bB_2^\top\\\hline
-\bC_1^\top\hat \bS^{-1}\bC_1&0&-(\bA_{11}-\bB_1\hat \bS^{-1}\bC_1)^\top&-(\bA_{21}-\bB_2\hat \bS^{-1}\bC_1)^\top\\
0&0&0&-\bA_{22}^\top\rix.
\end{eqnarray*}
The spectrum of $\wt\bH_{11}$ is the union of the spectra of
the submatrices $\bH_{11}$,
$\bA_{22}$, and $-\bA_{22}$. Since $\bA_{22}$ is asymptotically stable, if
$\bH_{11}$
does not have purely imaginary eigenvalues, so is
$\wt \bH_{11}$. Hence, $\wt\bH_{11}$ has a Lagrangian invariant subspace.
Consider the Riccati equation
\eq{pertARE}
\wt\Psi_{11} +\widetilde{\boldsymbol{\Xi}}_{11}=0,
\en
with $\widetilde{\boldsymbol{\Xi}}_{11}\ge 0$ being chosen such that $(\wt \bA_{11}-\wt \bB_1\hat \bS^{-1}\wt\bC_1,\wt\bC_1^\top\hat \bS^{-1}\wt\bC_1+\widetilde{\boldsymbol{\Xi}}_{11})$
is observable; it is clear that such a $\widetilde{\boldsymbol{\Xi}}_{11}$ always exists.
Recall that $(\wt \bA_{11}-\wt \bB_1\hat \bS^{-1}\wt\bC_1,\wt\bB_1)$
is controllable. The Riccati equation (\ref{pertARE}) corresponds to the Hamiltonian matrix
\[
\wt\bH_{11}(\widetilde{\boldsymbol{\Xi}}_{11}):=\wt \bH_{11}-\mat{c|c}0&0\\\hline \widetilde{\boldsymbol{\Xi}}_{11}&0\rix.
\]
For a sufficiently small (in norm) $\widetilde{\boldsymbol{\Xi}}_{11}$, by continuity, $\wt \bH_{11}(\widetilde{\boldsymbol{\Xi}}_{11})$ has a Lagrangian invariant subspace (e.g., when $\widetilde{\boldsymbol{\Xi}}_{11}$ is chosen small enough so that no eigenvalues of $\wt \bH_{11}(\widetilde{\boldsymbol{\Xi}}_{11})$ are  on the imaginary axis),
and by Lemma~\ref{lem:COSol}, $\wt \Psi=-\widetilde{\boldsymbol{\Xi}}_{11}\le 0$
has a positive definite solution $\wt\bQ_{11}$, where  $\wt\bQ_{11}$ can be chosen
so that all the eigenvalues of $\wt\bA_{11}-\wt\bB_1\bS^{-1}(\wt\bC_1-\wt\bB_1^\top\wt\bQ_{11})$ are in the closed left half complex plane.

If $\bH_{11}$ as given in (\ref{hatH}) has purely imaginary eigenvalues (potentially including $0$)
and we assume that it has a Lagrangian invariant subspace, then (\ref{RIE11}), written with equality,
still has a solution $Q_{11}^0>0$ such that all eigenvalues of
$\bA_{11}-\bB_1\hat \bS^{-1}(\bC_1-\bB_1^\top\bQ_{11}^0)$ are in the closed
left half complex plane.
Let $\wt \bQ_{11}^0=\mat{cc}\bQ_{11}^0&0\\0&0\rix\ge 0$. Then
$\wt\Psi_{11}(\wt\bQ_{11}^0)=0.$
Subtracting this
%
%
from (\ref{pertARE}) yields the Riccati equation
\eq{0ric}
(\wt\bA_{11}^0)^\top\wt\bY+\wt\bY\wt\bA_{11}^0
+\wt\bY\wt\bB_1\hat \bS^{-1}\wt\bB_1^\top\wt\bY+\widetilde{\boldsymbol{\Xi}}_{11}=0,
\en
where $\wt\bY=\wt\bQ_{11}-\wt\bQ_{11}^0$ and
\[
\wt\bA_{11}^0=\wt\bA_{11}-\wt\bB_1\hat \bS^{-1}
(\wt\bC_1-\wt\bB_1^\top\wt\bQ_{11}^0)
=\mat{cc}\bA_{11}-\bB_1\hat \bS^{-1}(\bC_1-\bB_1^\top\bQ_{11}^0)&0\\
\bA_{21}-\bB_2\hat \bS^{-1}(\bC_1-\bB_1^\top\bQ_{11}^0)&\bA_{22}\rix.
\]
Since the eigenvalues of $\bA_{11}-\bB_1\hat \bS^{-1}(\bC_1-\bB_1^\top\bQ_{11}^0)$
are in the closed half complex plane, there is
an invertible block lower triangular matrix
\[
\wt\bL=\mat{cc}\bL_{11}&0\\\bL_{21}&\bI\rix
\]
such that
\[
\wt\bL\wt\bA_{11}^0
\wt\bL^{-1}
=\wt\bL\mat{cc}\bA_{11}-\bB_1\hat \bS^{-1}(\bC_1-\bB_1^\top\bQ_{11}^0)
&0\\\bA_{21}-\bB_2\hat \bS^{-1}(\bC_1-\bB_1^\top\bQ_{11}^0)&\bA_{22}\rix\wt\bL^{-1}
=\mat{c|cc}\boldsymbol{\Sigma}_1&0&0\\\hline 0&\boldsymbol{\Sigma}_2&0\\0&\boldsymbol{\Sigma}_{32}&\bA_{22}\rix
=:\mat{cc}\boldsymbol{\Sigma}_1&0\\0&\boldsymbol{\Sigma}_2^0\rix,
\]
where $\boldsymbol{\Sigma}_1$ has only purely imaginary eigenvalues and
$\boldsymbol{\Sigma}_2^0$ is asymptotically stable.
The block $\bL_{11}$ is used for the similarity transformation
\[
\bL_{11}(\bA_{11}-\bB_1\hat \bS^{-1}(\bC_1-\bB_1^\top\bQ_{11}^0)\bL^{-1}_{11}
=\mat{cc}\boldsymbol{\Sigma}_1&0\\0&\boldsymbol{\Sigma}_2\rix,
\]
while $\bL_{21}$ is used to uncouple the diagonal blocks $\boldsymbol{\Sigma}_1$ and $\boldsymbol{\Sigma}_2^0$.
Repartition
\[
\wt\bL\wt \bB_1=\mat{c}\bB_1^0\\\bB_2^0\rix
\]
according to $\wt\bL\wt\bA_{11}^0\wt\bL^{-1}$.
We look for a solution of (\ref{0ric}) of the form
\[
\wt\bL^{-\top}\wt\bY\wt\bL^{-1}=\mat{cc}0&0\\0&\bY_2\rix.
\]
By taking a congruence transformation on both sides of (\ref{0ric}) with
$\wt\bL^{-\top}$ on the left and $\wt\bL^{-1}$ on the right, the resulting equation
reduces to
\[
(\boldsymbol{\Sigma}_2^0)^\top\bY_2+\bY_2\boldsymbol{\Sigma}_2^0+\bY_2\bB_2^0\hat \bS^{-1}(\bB_2^0)^\top\bY_2
+\boldsymbol{\Xi}_{22}=0
\]
for a suitably chosen $\boldsymbol{\Xi}_{22}\ge 0$, and $\widetilde{\boldsymbol{\Xi}}_{11}=\wt \bL^\top\mat{cc}0&0\\0&\boldsymbol{\Xi}_{22}\rix\wt\bL\ge 0$.
Since $(\wt\bA_{11}^0,\wt\bB_1)$ is controllable, so is
$(\boldsymbol{\Sigma}_2^0,\bB_2^0)$. Recall also that $\boldsymbol{\Sigma}_2^0$ is asymptotically stable.
Analogous to the previous case one can
choose (a sufficiently small) $\boldsymbol{\Xi}_{22}\ge 0$ with $(\boldsymbol{\Sigma}_2^0,\boldsymbol{\Xi}_{22})$ observable
and then the reduced Riccati equation has a
positive definite solution $\bY_2$ with the eigenvalues of
$\boldsymbol{\Sigma}_2^0+\bB_2^0\hat \bS^{-1}(\bB_2^0)^\top\bY_2$ in the
closed left half complex plane. Then
\[
\wt \bQ_{11}=\wt\bQ_{11}^0+\wt \bL^\top\mat{cc}0&0\\0&\bY_2\rix\wt\bL
\]
solves (\ref{pertARE}). Since
\begin{eqnarray*}
\wt\bL(\wt\bA-\wt \bB_1\hat \bS^{-1}(\wt\bC_1-\wt\bB_1^\top\wt \bQ_{11}))\wt\bL^{-1}
&=&\wt\bL\wt\bA_{11}^0\wt\bL+\wt\bL\wt\bB_1\hat \bS^{-1}\wt(\bL\wt\bB_1)^\top\mat{cc}0&0\\0&\bY_2\rix\\
&=&\mat{cc}\boldsymbol{\Sigma}_1&\bB_1^0\hat \bS^{-1}(\bB_2^0)^\top\bY_2\\
0&\boldsymbol{\Sigma}_2^0+\bB_2^0\hat \bS^{-1}(\bB_2^0)^\top\bY_2\rix,
\end{eqnarray*}
the eigenvalues of $\wt\bA-\wt \bB_1\hat \bS^{-1}(\wt\bC_1-\wt\bB_1^\top\wt \bQ_{11})$ are
in the closed left half complex plane. Because $\bQ_{11}^0>0$, $\bY_2>0$ and
\begin{eqnarray*}
\wt\bL^{-\top}\wt \bQ_{11}\wt\bL^{-1}&=&\wt\bL^{-\top}\wt\bQ_{11}^0\wt\bL^{-1}+\mat{cc}0&0\\0&\bY_2\rix
=\mat{cc}\bL_{11}^{-\top}\bQ_{11}^0\bL_{11}^{-1}&0\\0&0\rix+\mat{cc}0&0\\0&\bY_2\rix\\
&=&\mat{ccc}\bQ_{11}^0&\bQ_{12}^0&0\\
(\bQ_{12}^0)^\top&\bQ_{22}^0&0\\0&0&0\rix
+\mat{ccc}0&0&0\\0&\widehat \bY_{11}&
\widehat \bY_{12}\\0&\widehat \bY_{12}^\top&\widehat \bY_{22}\rix,
\end{eqnarray*}
one has $\wt\bQ_{11}>0$. Hence, in either case the equation (\ref{pertARE}) has
a positive definite solution $\wt\bQ_{11}$ with all the eigenvalues of
$\wt\bA-\wt \bB_1\hat \bS^{-1}(\wt\bC_1-\wt\bB_1^\top\wt \bQ_{11})$ in the closed
half complex plane for some
$\widetilde{\boldsymbol{\Xi}}_{11}\ge 0$.

Once we have such a solution $\wt\bQ_{11}$, we can solve
$\wt\Psi_{12}=0$ for $\wt\bQ_{12}$. This Sylvester equation
has a unique solution because $\wt\bA_{22}=\bA_{33}$ is asymptotically stable
and the eigenvalues  of $\wt\bA_{11}-\wt\bB_1\bS^{-1}(\wt\bC_1-\wt\bB_1^\top
\wt\bQ_{11})$ are in the closed left half complex plane.

Having solved the equality $\wt\Psi_{12}=0$ (but omitting details), we finally approach the solution of
the inequality $\wt\Psi_{22}\le 0$. We may consider the Lyapunov equation
\[
\wt\Psi_{22}=-\widetilde{\boldsymbol{\Xi}}_{22}\le 0.
\]
Since $\wt A_{22}=A_{33}$ is asymptotically stable, for any $\widetilde{\boldsymbol{\Xi}}_{22}\ge 0$ it
always has a solution $\wt\bQ_{22}$.
Then $\bQ_1=\bU\mat{cc}\wt \bQ_{11}&\wt\bQ_{12}\\\wt\bQ_{12}^\top&\wt\bQ_{22}\rix\bU^\top$
solves the Riccati equation
\[
\Psi(\bQ_1)=-\bU\mat{cc}\widetilde{\boldsymbol{\Xi}}_{11}&0\\0&\widetilde{\boldsymbol{\Xi}}_{22}\rix\bU^\top\le 0.
\]
With the assumption that $\bA_1-\bB_1\hat \bS^{-1}\bC_1$  is asymptotically stable,
this $\bQ_1$ must be positive semidefinite. Suppose that $\bQ_1\bU\bx=0$ for some $\bx=\mat{c}\bx_1\\\bx_2\rix\ne 0$. The vector $\bx_2\ne 0$, since otherwise $\wt\bQ_{11}\bx_1=0$ contradicting
the positive definiteness of $\wt\bQ_{11}$.
From $\Psi(\bQ_1)=-\bU\mat{cc}\widetilde{\boldsymbol{\Xi}}_{11}&0\\0&\widetilde{\boldsymbol{\Xi}}_{22}\rix\bU^\top$, one has
\[
\bx^\top\mat{cc}\widetilde{\boldsymbol{\Xi}}_{11}&0\\0&\widetilde{\boldsymbol{\Xi}}_{22}\rix\bx=\bx_1^\top\widetilde{\boldsymbol{\Xi}}_{11}\bx_1+\bx_2^\top\widetilde{\boldsymbol{\Xi}}_{22}\bx_2=0.
\]
But if we choose $\widetilde{\boldsymbol{\Xi}}_{22}>0$, this is not possible and thus $\bQ_1$ must be positive definite. Therefore,
choosing a positive definite $\widetilde{\boldsymbol{\Xi}}_{22}$ guarantees the corresponding solution $\bQ_1$ to be positive definite.

\begin{remark}\label{rem:Xi12}
{\rm
In the described construction of positive definite solutions to \eqref{LMIcond}, we have chosen
$\widetilde{\boldsymbol{\Xi}}_{11}$ ($\boldsymbol{\Xi}_{22}$) and $\widetilde{\boldsymbol{\Xi}}_{22}$
in order to guarantee the existence of positive definite solutions of the individual occurring Riccati and Lyapunov equations.

In a more general framework we can choose such a perturbation to turn the inequality in \eqref{LMIcond} into an equality by adding a positive semidefinite matrix to the left hand side to guarantee the existence of positive definite solutions. The set of all positive semidefinite perturbations of this kind will then help to characterize the solution set of  \eqref{LMIcond}.  To do this in detail is beyond the scope of this paper and has recently been investigated in \cite{MehX24}.
}
\end{remark}
We summarize the conditions for the existence of a positive definite
solution of
the matrix inequality (\ref{LMIred}) in the following theorem.
\begin{theorem}\label{lmisol_thm} Consider a  general  system of the form (\ref{GenSys}) with $ \bA$ stable and
$\hat \bS=\bD+\bD^\top>0$. Let $\bM$ be invertible such that
\[
\bM\bA\bM^{-1}=\mat{cc}\bA_1&0\\0&\bA_2\rix,\quad
\bM\bB=\mat{c}\bB_1\\\bB_2\rix,\quad\bC\bM^{-1}
=\mat{cc}\bC_1&\bC_2\rix,
\]
where $\bA_2$ is diagonalizable and contains all the purely imaginary
eigenvalues of $\bA$. Let $\bA_1,\bB_1,\bC_1$ have the
condensed form (\ref{defa1}) with an orthogonal matrix $\bU$. Then the matrix inequality (\ref{LMIred}) has a positive
definite solution $\bQ$ if and only if the following conditions hold.
\begin{enumerate}
\item[(a)]  There exists a positive definite matrix $\bQ_2$ satisfying
\[
\bB_2^\top\bQ_2=\bC_2,\quad \bA_2\bQ_2=\bQ_2\bA_2.
\]
\item[(b)] The block $\bA_{11}-\bB_1\hat \bS^{-1}\bC_1$ is asymptotically
stable.
\item[(c)] The Hamiltonian matrix $\bH_{11}$ defined in (\ref{hatH}) has
a Lagrangian invariant subspace.
\end{enumerate}

If these conditions hold, then the linear matrix inequality (\ref{LMIred}) has a positive definite solution
of the form $\bQ=\bM^{\top}\mat{cc}\bQ_1&0\\0&\bQ_2\rix \bM>0$,
where $\bQ_1>0$ solves (\ref{RIEred}),
$\bQ_2$ is determined from condition (a),
and
\begin{eqnarray*}
&&(\bA-\bB\hat \bS^{-1}\bC)^\top\bQ
+\bQ(\bA-\bB\hat \bS^{-1}\bC)+\bQ\bB\hat \bS^{-1}\bB^\top\bQ
+\bC^\top\hat \bS^{-1}\bC\\
&&=-\bM^\top\mat{cc}\bU\mat{cc}\widetilde{\boldsymbol{\Xi}}_{11}&0\\0&\widetilde{\boldsymbol{\Xi}}_{22}\rix\bU^\top&0\\0&0\rix\bM
\le 0,
\end{eqnarray*}
where $\widetilde{\boldsymbol{\Xi}}_{11}, \widetilde{\boldsymbol{\Xi}}_{22}$ are chosen as in the above explicit construction.
\end{theorem}
\proof
The proof follows from the explicit construction.
\eproof
\begin{remark}\label{rem_even} {\rm
If $\hat \bS>0$ then the eigenvalues of the Hamiltonian matrix ${\bf H}$ in (\ref{hampas}) are the finite eigenvalues  of
the \emph{regular even matrix pencil}
\begin{equation}\label{evenpencildef}
\lambda \boldsymbol{\mathcal{N}}-\boldsymbol{\mathcal{M}}:= \lambda \mat{ccc} 0 & \bI & 0 \\ -\bI & 0 & 0 \\ 0 & 0 & 0\rix-
\mat{ccc} 0 &  \bA & \bB \\  \bA^\top & 0 & \bC^\top \\ \bB^\top & \bC &  \hat \bS\rix
\end{equation}
which has \emph{index at most one}, i.e., the eigenvalues at $\infty$ are all semi-simple. It is straightforward to verify
that the relation (\ref{haminv}) is equivalent to
\begin{equation}\label{evenpencil}
\boldsymbol{\mathcal{N}}
\mat{c} -\bW_2 \\ \bW_1 \\ \bY \rix
\bZ=
\boldsymbol{\mathcal{M}}
\mat{c} -\bW_2\\ \bW_1 \\ \bY \rix
\end{equation}
with $\bY=\hat \bS^{-1}(\bB^\top\bW_2-\bC\bW_1)$, i.e.,
the columns of $\mat{ccc}-\bW_2^\top&\bW_1^\top&\bY^\top\rix^\top$ span the $n$-dimensional
deflating subspace of $\lambda\boldsymbol{\mathcal{N}}-\boldsymbol{\mathcal{M}}$ corresponding to the finite eigenvalues contained in $\bZ$.

A pencil $\lambda\, \boldsymbol{\mathcal{N}}-\boldsymbol{\mathcal{M}}$ is called \emph{even} if $\boldsymbol{\mathcal{N}}=-\boldsymbol{\mathcal{N}}^\top$ and
$\boldsymbol{\mathcal{M}}=\boldsymbol{\mathcal{M}}^\top$.
 Since even pencils
have the Hamiltonian spectral symmetry in the finite eigenvalues, see \cite{ByeMX07}, this means that there are equally many
eigenvalues in the open left and  in the open right half plane.

Numerically, to compute $\bQ(=\bW_2\bW_1^{-1})$ it is preferable to work with the even pencil (\ref{evenpencildef})  rather than with the Hamiltonian
matrix (\ref{hampas}), since explicit inversion of $ \hat \bS$ is avoided. Numerically stable structure preserving methods for computing an orthonormal matrix $\mat{ccc}-\bW_2^\top&\bW_1^\top&\bY^\top\rix^\top$
are available, see \cite{BenBMX99,BenLMV15,BruM08}.
}
\end{remark}
\begin{remark}\label{rem:passalone}{\rm
To show that  the system (\ref{GenSys}) is passive, it is sufficient that the linear matrix inequality (\ref{LMIred}) has a positive semidefinite solution $\bQ$.
In this case the conditions for the existence of solutions to~\eqref{LMIcond} can be relaxed. First of all, the
condition $\bB_2^\top\bQ_2=\bC_2$
can be relaxed to
$\mathsf{Ker }\,\bB_2\subseteq \mathsf{Ker }\,\bC_2^\top$,
$\rank\bC_2\bB_2=\rank\bC_2$, and $\bC_2\bB_2\ge 0$.
Also, $\bA_2$ may  have purely
imaginary eigenvalues with Jordan blocks. For example in the extreme case when $\bC_2=0$,
$\bQ_2=0$ always satisfies the conditions for any $\bA_2$.

Secondly, for (\ref{LMIredred}) or (\ref{RIEred}),
we still require that $\bA_{11}-\bB_1\hat \bS^{-1}\bC_1$ is asymptotically stable and
that the Hamiltonian matrix $\bH_{11}$ in (\ref{hatH}) has a Lagrangian invariant subspace.
Since this only requires that $\bQ_1 \ge 0$, a solution can be
determined in a simpler way.
We may simply set
$\wt \bQ_{11}=\wt\bQ_{11}^0$ and $\wt\Psi_{11}=0$.
Then with the block structure the solution of $\wt\Psi_{12}=0$
has a form $\wt\bQ_{12}=\mat{c} \bQ_{12}\\0\rix$. To solve
$\wt\Psi_{22}=0$ for $\wt\bQ_{22}$, one can show
$\bQ_1=\bU\mat{cc}\wt\bQ_{11}&\wt\bQ_{12}\\\wt\bQ_{12}^\top&
\wt\bQ_{22}\rix\bU^\top\ge 0$ and solves the Riccati equation
$\Psi(\bQ_1)=0$.
Also, in some circumstances the condition that $\wt\bA_{22}=\bA_{33}$
is asymptotically stable can be relaxed. In these relaxed cases however, it is necessary to transform the system to a descriptor formulation, see \cite{CheGH22}.
}
\end{remark}
{
 In this section we have recalled and extended classical results of the solution of the KYP inequality and associated Riccati inequality and have discussed the solution in the case that the system is not minimal but the symmetric part of the feedthrough term is positive definite. In the next section we discuss the general case and  we use transformations that can be implemented as numerically reliable procedures.
\section{Numerical procedures for the construction of port-Hamiltonian realizations}\label{sec:suff}
}
In the last section we have seen that the existence of a port-Hamiltonian realization for (\ref{GenSys})
reduces to the existence of a nonsingular matrix
$\bT\in \bRnn$ or a positive definite matrix $\bQ=\bT^\top\bT$ such that the matrix inequality 
(\ref{LMIcond}) holds. Note that since $\bQ=\bT^\top\bT$,
(\ref{LMIcond}) is equivalent to
\eq{sumineq}
\bW+\bW^\top\ge 0,
\en
where
\[
\bW=\mat{cc}\bT&0\\0&\bI\rix\mat{rr}-\bA&-\bB\\\bC&\bD\rix\mat{cc}\bT&0\\0&\bI\rix^{-1},
\]
for an invertible matrix $\bT$.

We develop here a constructive procedure to check these conditions. In our previous considerations the matrix $\bV$ that is used for a basis change in the input space  need only be invertible,  but to implement the transformations in a numerical stable manner,  we will require in the following that $\bV$ be real orthogonal.

{The procedures presented in the following three subsections follow from the procedures in \cite{WeiWS94}. The only extra result is an explicit formulation of the part of the matrix $\bT$ for dealing with the singularity of $\bD+\bD^\top$. We provide proofs for completeness.}
\subsection{The case that $\hat \bS= \bD+\bD^\top\geq 0$ is singular.}\label{singularS}
Suppose that the matrix
$\hat \bS= \bD+\bD^\top\geq 0$ is singular.
Consider an orthogonal matrix $\bV_0=[\bV_{0,1}, \bV_{0,2}]$, where $\bV_{0,1}$ is chosen so that
its columns form an orthonormal basis of the kernel of $ \hat \bS$.
To construct such a $\bV_0$ we can use a singular value or rank-revealing $QR$ decomposition, \cite{GolV96}.
Then
we have
\begin{equation}\label{feedthroughreduction}
\hat \bS=\bV_0^\top \bD\bV_0+\bV_0^{\top} \bD^\top\bV_0=
\bV^\top_0\left( \bD+ \bD^\top\right)\bV_0=\mat{cc} 0 & 0 \\ 0 & \bS_2\rix,
\end{equation}
where {\color{black}$0<\bS_2=\bD_2+\bD_2^\top$,
$\bD_2=\bV_{0,2}^\top \bD\bV_{0,2}\in\mathbb{R}^{s\times s}$
and $s=\rank(\bD+\bD^\top)$}.
Set
\begin{equation} \label{bcsplit}
\mat{cc} \bB_1 &  \bB_2 \rix = \bB \bV_0,\ \mat{cc}  \bC_1^\top &  \bC_2^\top \rix :=  \bC^\top \bV_0,
\end{equation}
each partitioned compatibly
with  {$\hat \bS$} as in (\ref{feedthroughreduction}).

Scaling the second block row and column of the matrix inequality (\ref{sumineq}) with $\bV_0^\top$ and $\bV_0$ respectively, we obtain the matrix inequality
\begin{equation} \label{LMItrans}
\left[\begin{array}{ccc} -(\bT\bA\bT^{-1})^{\top} -\bT\bA\bT^{-1} & -\bT\bB_1 +(\bC_1\bT^{-1})^\top &-\bT\bB_2 +(\bC_2\bT^{-1})^\top\\
-(\bT\bB_1)^\top +\bC_1\bT^{-1} & 0 & 0\\
-(\bT\bB_2)^\top+\bC_2\bT^{-1} & 0& \bS_2 \\
\end{array}\right]\ge 0
\end{equation}
which has an invertible solution $\bT$ if and only if  the matrix inequality
\begin{equation} \label{LMIconst}
\left[\begin{array}{cc} -(\bT\bA\bT^{-1})^{\top} -\bT\bA\bT^{-1} &-\bT\bB_2 +(\bC_2\bT^{-1})^\top\\
-(\bT\bB_2)^\top+\bC_2\bT^{-1} & \bS_2 \\
\end{array}\right]\geq 0
\end{equation}
has an invertible solution $\bT$ satisfying the constraint $(\bT\bB_1)^\top -\bC_1\bT^{-1}=0$.  We characterize conditions when this constraint is satisfied in the following subsections.
Notice that we are, in effect, restricting the input and output space to the invertible part of $ \bD+ \bD^\top$.  Once these restricted transformation matrices have been constructed, full transformations  satisfying the given constraint can be obtained by extending to the full space.
{\color{black}
\begin{remark}\label{rem:almostsingular}{\rm If $\hat \bS$ is positive definite but nearly singular, then the coefficients of the Riccati inequality in \eqref{RIE} may incur large relative errors when formed. In this case, it may be appropriate to regularize the problem by perturbing $\hat \bS$ to a nearby positive semidefinite (but now singular) problem as in \eqref{feedthroughreduction}, where the matrix $\bS_2$ is now well conditioned with respect to inversion.
In our procedure we do this by setting small positive eigenvalues of $\hat \bS$ to zero.
}
\end{remark}
}
\subsection{Construction  in the case
$ \bD=-  \bD^\top$}
To explicitly construct the transformation to port-Hamiltonian form  let us first discuss the extreme case that $\hat \bS=0$, i.e., that $ \bD= - \bD^\top$. Considering the matrix $\bK$ in \eqref{Kdef}, to satisfy $\bK\geq 0$, we must have $\bP=0$, and the block $\bV_{0,2}$ in the transformation of the feedthrough term is void, while $\bV_0=\bV_{0,1}$ is any orthogonal matrix. 
%
\begin{corollary}\label{Prop0case}
For a state-space system of the form \eqref{GenSys} with
$ \bD+  \bD^\top=0$ the following two statements are equivalent:
\begin{enumerate}
\item  There exists a change of basis $\bx= \bT^{-1}\bfxi$ with an invertible matrix $\bT\in \bRnn$ such that the resulting realization in the new basis {together with the same $\bu$ and $\by$ ($\bV=\bI$ in
(\ref{IOtransformations}))} has port-Hamiltonian structure as in (\ref{FullpHSys}).
\item There exists an invertible matrix $\bT$ such that
\begin{equation} \label{pHSuffCondDeq0}
a)\ (\bT{\bB})^\top ={\bC}\bT^{-1} \qquad \mbox{and} \qquad b)\ (\bT{\bA}\bT^{-1})^\top+\bT{\bA}\bT^{-1}\leq 0.
\end{equation}
\end{enumerate}
\end{corollary}
\proof {The proof follows directly from \eqref{LMIconst} and \eqref{LMItrans} and the preceding results.}
\eproof

Note that without the constraint (\ref{pHSuffCondDeq0}a),
if $ \bA$ is stable, then by Lemma~\ref{LyapStab},
the second condition (\ref{pHSuffCondDeq0}b) can always  be satisfied.
Adding the constraint (\ref{pHSuffCondDeq0}a), however, makes the question nontrivial.

We have the following characterization of the transformation matrices $\bT$
that satisfy (\ref{pHSuffCondDeq0}a).
\begin{lemma}{\cite{WeiWS94}}\label{lemT} Consider $ \bB,\bC^\top\in {\mathbb R}^{n\times m}$,
and assume that $\rank \bB=r$.
\begin{enumerate}
\item[a)] There exists an invertible transformation $\bT$ satisfying condition (\ref{pHSuffCondDeq0}a)
if and only if $\mathsf{Ker } \,\bC^\top =\mathsf{Ker }\,\bB$,
$\rank\bC\bB=r$ and $\bC\bB\ge 0$, or equivalently, there exists
an invertible (orthogonal) matrix $ \bW$ such that
\[
\bB \bW=\mat{cc}\bB_1&0\rix,\quad
\bC^\top\bW=\mat{cc}\bC_1^\top&0\rix,
\quad \bC_1\bB_1= \bY \bY^\top>0,
\]
where $\bB_1,\bC_1^\top\in{\mathbb R}^{n\times r}$ have full column rank
and $\bY\in{\mathbb R}^{r\times r}$ is invertible.
\item[b)] Let $\bN_{\bB}\in \mathbb R^{n\times (n-r)}$ have columns that form a basis of
$\mathsf{Ker }\,\bB^\top$. If Condition a) is satisfied,
then any $\bT$ satisfying condition (\ref{pHSuffCondDeq0}a) has the form
$\bT= \bU \bT_{\bZ} \bT_0$ with
\begin{equation} \label{Tform}
 \bT_0=\left[\begin{array}{c} \bN_{ \bB}^\top\\
 \bY^{-1}\bC_1 \end{array}\right],\
 \bT_{ \bZ}=\left[\begin{array}{cc} \bZ & 0 \\ 0 & \bI \end{array}\right],
\end{equation}
where $ \bU\in {\mathbb R}^{n\times n}$ is an arbitrary orthogonal
matrix and $ \bZ\in {\mathbb R}^{(n-r)\times (n-r)}$ is an arbitrary nonsingular matrix.
\end{enumerate}
\end{lemma}
\proof
Condition (\ref{pHSuffCondDeq0}a) is equivalent to
$ \bC=  \bB^\top \bT^\top \bT$ and
the following conditions will be necessary for the existence of an invertible $\bT$ with this property:
\begin{eqnarray}
\mathsf{Ker }\, \bC^\top&=&\mathsf{Ker }\,\bT^\top\bT \bB
=\mathsf{Ker}\,  \bB,\nonumber \\
\label{nscond1}
0&\leq &\bC\bB=\bB^\top \bT^\top \bT \bB,\\
\rank \bC\bB &=&\rank \bB^\top \bT^\top \bT \bB=\rank\bB=r.
\nonumber
\end{eqnarray}

Conversely,  observe that $\mathsf{Ker}\,\bC^\top=\mathsf{Ker}\, \bB$ is equivalent
to the existence of an orthogonal matrix $\bW\in \mathbb R^{m\times m}$ such that
\begin{equation}\label{bc1}
\bB\bW=\mat{cc}\bB_1&0\rix,\quad
\bC^\top\bW=\mat{cc}\bC_1^\top&0\rix,
\end{equation}
with $\bB_1,\bC_1^\top\in{\mathbb R}^{n\times r}$ having full column rank.
The conditions $\rank\bC\bB=r$ and $\bC\bB\ge 0$ together
are equivalent to {\color{black}$\bC_1\bB_1>0$}.
Thus, there must exist an invertible matrix $\bY$ (e.g., a Cholesky factor or a positive-definite square root) such that
$\bC_1 \bB_1=\bY\bY^\top$.

The matrix $ \bT_0$ as in (\ref{Tform}) is then well defined.
Furthermore, $ \bT_0$ is invertible, since if $ \bT_0\by=0$ for some
vector $\by$, then $\bC_1\by=0$ and
$\bN_{\bB}^\top\by=0$.  The latter statement implies that
$\by\in \mathsf{Ran}(\bB)$, so $\by=\bB_1 \bz$ for some $\bz$
and, furthermore, $\bC_1\bB_1\bz=0$.
This in turn implies that $\bz=0$ and $\by=0$; so $ \bT_0$ is injective, and hence invertible.

The invertibility of $ \bT_0$ implies that
\begin{eqnarray*}
 \bB^\top  \bT_0^\top  \bT_0 &=& \bW^{-\top}\mat{c} \bB_1^\top\\0\rix
 \left[\bN_{\bB}\quad \bC_1^\top  \bY^{-\top} \right]
  \left[\begin{array}{c}  \bN_{\bB}^\top\\  \bY^{-1}\bC_1  \end{array}\right]\\
  &=&  \bW^{-\top}\mat{cc}0&( \bY^{-1}\bC_1\bB_1)^\top \\  0&0 \rix
  \left[\begin{array}{c} \bN_{\bB}^\top \\  \bY^{-1}\bC_1 \end{array}\right]\\
   &=&  \bW^{-\top}\mat{c}(\bC_1\bB_1)(\bC_1\bB_1)^{-1}\bC_1\\0\rix
   = \bW^{-\top}\mat{c}\bC_1\\0\rix =\bC.
\end{eqnarray*}
Hence, (\ref{pHSuffCondDeq0}a) holds with $\bT= \bT_0$.

Now suppose that $\bT$ is any invertible transformation satisfying (\ref{pHSuffCondDeq0}a).  Then,
\[
\bB^\top \bT^\top (\bT  \bT_0^{-1}) =\bC \bT_0^{-1}=\bB^\top\bT_0^\top,
\]
which is equivalent  to
\[
\bB^\top_1  \bT_0^\top\left( (\bT  \bT_0^{-1})^\top(\bT  \bT_0^{-1})-\bI \right) =
\left[ 0\quad\bY  \right]\left( (\bT  \bT_0^{-1})^\top(\bT  \bT_0^{-1})-\bI \right)=0.
\]
 From this, it follows that
\[
(\bT  \bT_0^{-1})^\top(\bT  \bT_0^{-1})
= \left[\begin{array}{cc} \bZ^\top\bZ & 0 \\ 0 & \bI \end{array}\right]
 \]
for some invertible matrix $ \bZ$,
and that $\bT  \bT_0^{-1}\left[\begin{array}{cc} \bZ^{-1} & 0 \\ 0 & \bI \end{array}\right]= \bU$ must be real orthogonal.
\eproof

In order to explicitly construct a transformation matrix $ \bT_0$ as in (\ref{Tform}), it will be useful to construct bi-orthogonal bases for the two subspaces $\mathsf{Ker }\,\bB^\top$
and $\mathsf{Ker}\,\bC$.  Toward this end, let $\bN_{\bC}$ contain columns that form a basis of $\mathsf{Ker}\,\bC$, so that $\mathsf{Ran}\,\bN_{\bC}=\mathsf{Ker}\, \bC$. Such a matrix is easily constructed in a numerically stable way via the singular value decomposition or a rank-revealing $QR$ decomposition of $ \bC$, see \cite{GolV96}. Since $ \bB$ and $ \bC^\top$ are assumed to satisfy $\mathsf{Ker}\,\bB=\mathsf{Ker}\,\bC^\top$,
we have singular value decompositions
\begin{equation}\label{SVDsBC}
 \bB =[\bU_{ \bB,1}\quad \bU_{\bB,2}] \mat{cc}\boldsymbol{\Sigma}_{\bB}&0 \\ 0&0
\rix {\bV_{ \bB}}^\top,\quad
 \bC =\bU_{\bC} \mat{cc}\boldsymbol{\Sigma}_{ \bC} & 0\\0&0 \rix
\mat{c}\bV_{\bC,1}^\top\\ \bV_{\bC,2}^\top \rix
\end{equation}
with $\boldsymbol{\Sigma}_{\bB},\boldsymbol{\Sigma}_{\bC}\in{\mathbb R}^{r\times r}$
both invertible.
We then obtain
\begin{equation} \label{NBNC}
\bN_{\bC}= \bV_{ \bC,2},\quad \bN_{\bB}= \bU_{ \bB,2}.
\end{equation}
Observe that $\bN_{\bB}^\top\bN_{\bC}$ is nonsingular,
since if $\bN_{\bB}^\top\bN_{\bC}\by=0$ and $\bz=\bN_{\bC}\by$, then $\bN_{\bB}^\top\bz=0$ implies that
$\bz\in \mathsf{Ran}\,\bB$. But then, $\bz=\bB_1\bw
=\bN_{\bC}\by$
 implies $\bC_1\bB_1 \bw=0$,
  where $\bB_1$ and $\bC_1$ are defined in (\ref{bc1}), and so, $\bw=0$ and hence $\bz=0$. Then, since
$\bN_{\bC}$ has full column rank, we have $\by=0$. Thus, $\bN_{\bB}^\top\bN_{\bC}$ is injective, hence invertible.

Performing another singular value decomposition,  $\bN_{\bB}^\top\bN_{\bC}=\widetilde{\bU}\widetilde{\Delta}\widetilde{\bV}^\top$, with $\wt \Delta$ positive diagonal, and $\wt \bU, \wt \bV$ real orthogonal, we can perform a change of basis
$\wt \bN_{\bB}=\bN_{\bB}\widetilde{\bU} \widetilde{\Delta}^{-1/2}$ and $\wt \bN_{\bC}=\bN_{\bC} \widetilde{\bV} \widetilde{\Delta}^{-1/2}$ and obtain that the columns of $\wt \bN_{\bB}$ form a basis for $\mathsf{Ker}\,\bB^\top$,
the columns of $\wt \bN_{\bC}$ form a basis for $\mathsf{Ker}\,\bC$ and these two bases are bi-orthogonal, i.e.,
$\wt \bN_{\bB}^\top \wt \bN_{\bC} = \bI$, and we have
\eq{invT0}
 \bT_0=\mat{c}\wt \bN_{\bB}^\top\\  \bY^{-1}\bC_1\rix,
\quad
 \bT_0^{-1}=\mat{cc} \wt \bN_{\bC}&\bB_1 \bY^{-\top}\rix.
\en
%
%
Note that $\bW, \bB_1,\bC_1$ in Lemma~\ref{lemT} can be determined by the SVDs
in (\ref{SVDsBC}).

Using the formula (\ref{Tform}), we can express the conditions for a transformation to pH form that we have obtained so far in a more concrete way.
\begin{corollary}\label{lem0case}
 Consider system (\ref{GenSys}) with $ \bD=- \bD^\top$ and $\rank\bB=r$.
Let the columns of $\wt \bN_{\bB}$ and $\wt \bN_{\bC}$ span the kernels of $\bB^\top$, $\bC$ and satisfy
$\wt \bN_{\bB}^\top\wt \bN_{\bC}=\bI$.
Then system  (\ref{GenSys}) is equivalent to a pH system  if and only if
\begin{enumerate}
\item $\mathsf{Ker}\,\bC^\top=\mathsf{Ker}\,\bB$,
$\rank \bC\bB=r$, $\bC\bB\ge 0$,  and
\item there exists an invertible matrix $\bZ$ such that
\begin{equation} \label{LMIform2}
\mat{cc}\bZ&0\\0&\bI\rix
 \bT_0 \bA \bT_0^{-1}
\mat{cc}\bZ^{-1}&0\\0&\bI\rix+\left(\mat{cc}\bZ&0\\0&\bI\rix
 \bT_0 \bA \bT_0^{-1}
\mat{cc}\bZ^{-1}&0\\0&\bI\rix\right )^\top\leq 0,
\end{equation}
and  $ \bT_0$, $ \bT_0^{-1}$ are defined in (\ref{invT0}).
 %
\end{enumerate}
\end{corollary}
\proof The condition follows from Corollary~\ref{Prop0case} and
the representation (\ref{Tform}) by setting $ \bU=\bI$ and
$ \bT_0$ as in (\ref{invT0}).
\eproof
{\color{black}
\begin{remark}\label{rem:nonuniqueness}{\rm It is clear that the
matrices $\bW$ and $\bY$ in Lemma~\ref{lemT} are not unique.
In fact any matrix $\bW=\mat{cc}\bW_{11}&0\\\bW_{21}&\bW_{22}\rix$
with invertible $\bW_{11}$ and $\bW_{22}$ can replace $\bW$.
As a consequence, $\bB_1$ and $\bC_1$ will be replaced by
$\bB_1\bW_{11}$ and $\bW_{11}^\top\bC_1$. Then $\bY$ can be replaced by $\bW_{11}^\top \bY\bU_\bY$ for any real orthogonal matrix $\bU_\bY$, and $\bY^{-1}\bC_1$, $\bB_1\bY^{-\top}$
become $\bU_\bY^\top\bY^{-1}\bC_1$ and $\bB_1\bY^{-\top}\bU_\bY$.
The matrices $\wt \bN_\bB$ and $\wt \bN_\bC$ are not unique either. They can be replaced by $\wt\bN_\bB\bN^\top$ and
$\wt\bN_\bC\bN^{-1}$ with an arbitrary invertible matrix $\bN$.
In the end, $\bT_0$ can be replaced by
$\mat{cc}\bN&0\\0&\bU_\bY^\top\rix\bT_0$. The orthogonal matrix
$\bU_\bY$ can be absorbed in the matrix $\bU$ in the formula for
$\bT$ in Lemma~\ref{lemT}. The matrix $\bZ$ in $\bT_\bZ$ needs to be replaced by $\bZ\bN^{-1}$. However, $\bZ\wt \bN_\bB^\top$ is independent of $\bN$.
Note that $\bU$ affects $\bT$
but not $\bQ=\bT^\top\bT$.
}
\end{remark}
}
\subsection{Construction in the case of general $\bD$}

For the case that $\bD$ is general we will present a recursive procedure which is analogous to the index reduction procedure for differential-algebraic equations in \cite{KunM24}. The first step is to perform the transformations
(\ref{feedthroughreduction}), (\ref{bcsplit}), and to obtain the following
characterization when a transformation to port-Hamiltonian form~\eqref{FullpHSys} exists.

\begin{lemma}\label{lemgcase}
 Consider system (\ref{GenSys}) transformed as in (\ref{feedthroughreduction}) and (\ref{bcsplit}).
Then the system is equivalent to a port-Hamiltonian system of the form ~\eqref{FullpHSys} if and only if {\color{black}the following two conditions hold}
\begin{enumerate}
\item $\mathsf{Ker}\,\bC_1^\top=\mathsf{Ker}\,\bB_1$,
$\rank \bC_1\bB_1=\rank\bB_1{\color{black}=r}$, and $\bC_1\bB_1\ge 0$, or equivalently,
there exists an invertible (orthogonal) matrix $\bW$ such that
\eq{factbc1}
\bB_1\bW=\mat{cc}\wh \bB_1&0\rix,\quad \bC_1^\top\bW=\mat{cc}\wh\bC_1^\top&0\rix,\qquad
\wh \bC_1\wh \bB_1=\bY\bY^\top>0,
\en
where $\wh \bB_1,\wh\bC_1^\top\in\mathbb R^{n\times r}$ have full column rank and
$\bY\in\mathbb R^{r\times r}$ is invertible, and
\item there exists an invertible matrix $\bZ$ such that
%
%
%
\eq{deftG}
\wt\bY+\wt\bY^\top\ge 0,
\en
with
\eq{defF}
\wt\bY:=\mat{cc|c}\bZ&0&0\\0&{\bI_{r}}&0\\\hline 0&0&{\bI_{s}}\rix
\mat{c|c}-\bT_0\bA\bT_0^{-1}&-\bT_0\bB_2\\\hline
\bC_2\bT_0^{-1}&
{\color{black}\bD_2}\rix\mat{cc|c}\bZ&0&0\\
0&{\bI_{r}}&0\\\hline 0&0&{\bI_{s}}\rix^{-1}
\en
and
\begin{equation}\label{t0def}
\bT_0=\mat{c}\wt \bN_{\bB_1}^\top\\\bY^{-1}\wh\bC_1\rix,\quad
\bT_0^{-1}=\mat{cc}\wt \bN_{\bC_1}&\wh\bB_1\bY^{-\top}\rix,
\end{equation}
and
the columns of full rank matrices $\wt \bN_{\bB_1}$ and $\wt \bN_{\bC_1}$
form the kernels of $\bB_1^\top$, $\bC_1$ respectively, {\color{black}satisfying
$\wt \bN_{\bB_1}^\top\wt \bN_{\bC_1}=\bI_{n-r}$}.
\end{enumerate}
\end{lemma}
\proof
{Transforming the system to the form (\ref{feedthroughreduction}) and (\ref{bcsplit}), the
inequality \eqref{LMItrans} can only hold if the second block column and row are zero. This gives the characterization in part 1., since in this case
condition (\ref{pHSuffCondDeq0}a)  has the form $(\bT\bB_1)^\top=\bC_1\bT^{-1}$. In this way the linear matrix inequality only has to be solved in the part where the feedthrough term has a positive definite symmetric part. The solvability for this part is then characterized
by applying the result of Lemma~\ref{lemT} to
$\bB_1$ and $\bC_1$ and one obtains $\bT$ as the one in Lemma~\ref{lemT} b) with $\bT_0$ of the form as (\ref{Tform}).
The result is then proved by applying this formula to (\ref{sumineq}).}
\eproof



{In order to determine $\bT$ we still need to find
$\bZ$. Using the partition in the transformation matrix in (\ref{defF}) we obtain three blocks in $\wt\bY$  and we
can repartition the middle factor as}
%
\begin{eqnarray}
\nonumber
\mat{cc}-\bT_0\bA\bT_0^{-1}&-\bT_0\bB_2\\
\bC_2\bT_0^{-1}&{\bD_2}\rix
&=&\mat{c|cc}-\wt \bN_{\bB_1}^\top\bA\wt\bN_{\bC_1}&
-\wt\bN_{\bB_1}^\top\bA\wh\bB_1\bY^{-\top}&-\wt\bN_{\bB_1}^\top\bB_2\\
\hline
-\bY^{-1}\wh\bC_1 \bA\wt \bN_{\bC_1}&-\bY^{-1}\wh\bC_1\bA\wh\bB_1\bY^{-\top}&
-\bY^{-1}\wh\bC_1\bB_2\\
\bC_2\wt\bN_{\bC_1}&\bC_2\wh\bB_1\bY^{-\top}&
{\bD_2}\rix\\
\label{redsys}
&=:&\mat{rr}-\wt\bA_1&-\wt \bB_1\\ \wt\bC_1&\wt\bD_1\rix.
\end{eqnarray}
Thus, we have
\[
\wt\bY{\color{black}(\bZ)}=\mat{cc}\bZ&0\\0&{\bI_{r+s}}\rix
\mat{rr}-\wt\bA_1&-\wt \bB_1\\\wt\bC_1&\wt\bD_1\rix
\mat{cc}\bZ&0\\0&{\bI_{r+s}}\rix^{-1},
\]
and condition (\ref{deftG}) is exactly the same as the condition (\ref{sumineq}), just replacing
$\bT$, $\bA,\bB,\bC,\bD$
with $\bZ$, $\wt\bA_1$, $\wt\bB_1$, $\wt\bC_1$, $\wt\bD_1$.
Hence, the existence of $\bZ$ can be checked again by using
Lemma~\ref{lemgcase}.

This implies that  the procedure of checking the existence of a
transformation from (\ref{GenSys}) to a pH system can be performed in a recursive way. One first performs the transformation (\ref{feedthroughreduction}) and checks whether a condition as in Part 1 of
Lemma~\ref{lemgcase} {holds. If this is not the case then the transformation to port-Hamiltonian form does not exist. If the condition holds, then  one checks} whether
$\wt \bD_1+\wt \bD_1^\top$ in (\ref{redsys}) is {\color{black}positive definite}. If it is and if
the associated matrix inequality does not have a positive definite solution, then a transformation to {port-Hamiltonian form} does not exist.
Otherwise, a 
transformation matrix $\bT$
can be constructed by computing $\bZ$ satisfying $\wt\bY(\bZ)+\wt\bY(\bZ)^\top\ge 0$
and the matrix {\color{black}$\bT_\bZ=\mat{cc}\bZ&0\\0&\bI\rix$} is formed accordingly. In the remaining case, i.e.,
{if } a condition as in Part 1 of Lemma~\ref{lemgcase} holds and $\wt\bD_1+\wt\bD_1^\top$ is singular, {then} the process is repeated {with the  reduced problem}.

To formalize the recursive procedure, let
\[
\bG_0=\mat{rr}- \bA&-\bB\\\bC&\bD\rix,
\]
%
and suppose that (\ref{pHSuffCondDeq0}a) is satisfied. Then form
\[
\bG_1:=\wt \bT_0
\wt \bV_0^\top\bG_0\wt\bV_0\wt\bT_0^{-1},
\]
where
\[
\wt \bT_0=\mat{cc}\bT_0&0\\0&{\color{black}\bI_m}\rix,\quad
\wt\bV_0=\mat{cc}{\color{black}\bI_n}&0\\0 &\bV_0\rix,
\]
$\bV_0=\mat{cc}\bV_{0,2}&\bV_{0,1}\rix$ is the matrix  in the decomposition (\ref{feedthroughreduction})
times a permutation that interchanges the last block columns
and $\bT_0$ is obtained from~(\ref{t0def}). In this way we obtain
\[
\bG_1= \mat{c|cc|c}-\wt \bN_{\bB_1}^\top \bA\wt \bN_{\bC_1}
&-\wt\bN_{\bB_1}^\top\bA\wh\bB_1\bY^{-\top}&-\wt\bN_{\bB_1}^\top\bB_2&0\\ \hline
-\bY^{-1}\wh\bC_1 \bA\wt \bN_{\bC_1}&-\bY^{-1}\wh\bC_1\bA\wh\bB_1\bY^{-\top}&
-\bY^{-1}\wh\bC_1\bB_2&-{\color{black}\wt\Gamma}\\
\bC_2\wt\bN_{\bC_1}&\bC_2\wh\bB_1\bY^{-\top}&
\bV_{0,2}^\top\bD\bV_{0,2}&\bV_{0,2}^\top\bD\bV_{0,1}\\\hline
0&{\color{black}\wt\Gamma^\top}&
\bV_{0,1}^\top\bD\bV_{0,2}&\bV_{0,1}^\top\bD\bV_{0,1}
\rix,
\]
where, by using (\ref{factbc1}),
\begin{eqnarray*}
\bC_1\bT_0^{-1}&=&\mat{cc}0&\bW^{-\top}\mat{c}\wh\bC_1\\0\rix\wh \bB_1
\bY^{-\top}\rix=\mat{cc}0&\bW^{-\top}\mat{c}\bY\\0\rix\rix=\mat{cc}0&{\color{black}\wt\Gamma^\top}\rix,
\\
\bT_0\bB_1&=&\mat{c}0\\\bY^{-1}\wh \bC_1\mat{cc}\wh\bB_1&0\rix\bW^{-1}\rix=
\mat{c}0\\\mat{cc}\bY^\top&0\rix \bW^{-1}\rix=\mat{c}0\\{\color{black}\wt\Gamma}\rix.
\end{eqnarray*}
By (\ref{feedthroughreduction}), we have
\[\bV_{0,1}^\top\bD\bV_{0,2}=-(\bV_{0,2}^\top\bD\bV_{0,1})^\top,\qquad
\bV_{0,1}^\top\bD\bV_{0,1}=-(\bV_{0,1}^\top\bD\bV_{0,1})^\top.
\]
%
%
%
So we can express $\bG_1$ as
\[
\bG_1=\mat{rr|r}-\wt\bA_1&-\wt\bB_1&0\\\wt\bC_1&\wt\bD_1&-\Gamma_1\\\hline
0&\Gamma_1^\top&\Phi_1\rix,
\]
{\color{black}where $\wt\bA_1\in\mathbb{R}^{n_1\times n_1}$,
$\wt\bB_1,\wt\bC_1^\top\in\mathbb{R}^{n_1\times m_1}$,
$\wt\bD_1\in\mathbb{R}^{m_1\times m_1}$,
$\Gamma_1\in\mathbb{R}^{m_1\times d_1}$,
$\Phi_{1}=-\Phi_1^\top
\in\mathbb{R}^{d_1\times d_1}$
with $n_1=n-r$, $m_1=r+s$, and $d_1=m-s$.
}
If $\mat{rr}-\wt\bA_1&-\wt\bB_1\\\wt\bC_1&\wt\bD_1\rix$ satisfies
condition (\ref{pHSuffCondDeq0}a), then in an analogous way we construct
\[
\wt \bT_1=\diag(\bT_1,{\color{black}\bI_{m_1},\bI_{d_1}}),\quad \wt\bV_1=\diag({\color{black}\bI_{n_1}},\bV_1,{\color{black}\bI_{d_1}})
\]
such that
\[
\bG_2=\wt\bT_1\wt\bV_1^\top\bG_1\wt\bV_1\wt\bT_1^{-1}=
\mat{rr|rr}-\wt\bA_2&-\wt\bB_2&0&0\\
\wt\bC_2&\wt\bD_2&-\Gamma_2&-\Gamma_{11}\\\hline
0&\Gamma_2^\top&\Phi_2&-\Gamma_{12}\\
0&\Gamma_{11}^\top&\Gamma_{12}^\top&\Phi_1\rix,
\]
{\color{black}where $\wt\bA_2\in\mathbb{R}^{n_2\times n_2}$,
$\wt\bB_2,\wt\bC_2^\top\in\mathbb{R}^{n_2\times m_2}$,
$\wt\bD_2\in\mathbb{R}^{m_2\times m_2}$,
$\Gamma_2\in\mathbb{R}^{m_2\times d_2}$,
$\Phi_{2}=-\Phi_2^\top
\in\mathbb{R}^{d_2\times d_2}$
with $n_2=n_1-r_1$, $m_2=r_1+s_1$, $d_2=m_1-s_1$,
and $r_1=\rank \wt\bB_1$, $s_1=\rank (\wt \bD_1+\wt\bD_1^\top)$.
We then proceed.}

 {The procedure stops in two cases, (a)
condition (\ref{pHSuffCondDeq0}a) is violated or
$\wt \bD_j+\wt \bD_j^\top$ is indefinite for some $j$ and no transformation to port-Hamiltonian form exists, or (b)
after $k$ steps one arrives at}
\[
\bG_k=\wt\bT_{k-1}\wt \bV_{k-1}^\top\ldots\wt\bT_0\wt\bV_0^\top\bG_0
\wt\bV_0\wt\bT_0^{-1}\ldots\wt\bV_{k-1}\wt\bT_{k-1}^{-1}=
\mat{cc|c}-\wt\bA_k&-\wt\bB_k&0\\
\wt\bC_k&\wt\bD_k&-\wt\Gamma_k\\\hline
0&\wt \Gamma_k^\top&\wt\Phi_k\rix,
\]
where $\wt\bD_k+\wt\bD_k^\top$ is positive definite, $\wt\Phi_k=-\wt\Phi_k^\top$,
\[
\wt \bT_j=\diag(\bT_j,{\color{black}\bI_{m_j}},\bI),\quad
\wt \bV_j=\diag({\color{black}\bI_{n_j}},\bV_j,\bI),
\]
and {\color{black}$n_j$,$m_j$ are the sizes of $\bT_j$, $\bV_j$, respectively,} for $j=0,\ldots,k-1$.
{Case (b) happens because $n_j$ is strictly decreasing
until $r_{j}$ is void or zero for some $j$. If $r_{j}$ is void, then
$\wt \bD_j+\wt \bD_j$ is positive definite already. If $r_j=0$,
then $n_{j+1}=n_j$ and $m_{j+1}=s_j$ indicating that
$\wt \bD_{j+1}+\wt\bD_{j+1}^\top$ is positive definite.}

If there exists an invertible matrix $\bT_k$ such that
\[
\mat{cc}\bT_k&0\\0&\bI\rix
\mat{cc}-\wt\bA_k&-\wt\bB_k\\
\wt\bC_k&\wt\bD_k\rix\mat{cc}\bT_k^{-1}&0\\0&\bI\rix
+\left(\mat{cc}\bT_k&0\\0&\bI\rix
\mat{cc}-\wt\bA_k&-\wt\bB_k\\
\wt\bC_k&\wt\bD_k\rix\mat{cc}\bT_k^{-1}&0\\0&\bI\rix\right)^\top\ge 0,
\]
see the solvability conditions in the previous section,
then with $\wt\bT_k=\diag(\bT_k,\bI,\bI)$, we have
\begin{equation}\label{kstepineq}
\wt\bT_k\bG_k\wt\bT_k^{-1}+(\wt\bT_k\bG_k\wt\bT_k^{-1})^\top\ge 0.
\end{equation}
Observe that for each $i,j$ with $j\geq i$ the matrices  $\wt \bV_i$ and $\wt \bV_i^\top$ commute with
$\wt\bT_j$ and $\wt\bT_j^{-1}$, and thus setting
\[
\wt\bT=\wt\bT_k\ldots\wt\bT_0,\quad \wt\bV=\wt\bV_0\ldots\wt\bV_{k-1},
\]
then
\[
\wt\bT_k\bG_k\wt\bT_k^{-1}=
\wt\bV^\top\wt\bT\bG_0\wt\bT^{-1}\wt\bV,
\]
and inequality (\ref{kstepineq}) implies that
\[
\wt\bT\bG_0\wt\bT^{-1}+(\wt\bT\bG_0\wt\bT^{-1})\ge  0.
\]
Then the desired transformation matrix $\bT$ is positioned in the
top diagonal block of $\wt \bT$, and the matrix $\bV$ is positioned in the bottom diagonal block of
$\wt\bV$.

\begin{remark}
\label{rem:recursive}{\rm The recursive procedure described above requires at each step the
computation of three singular value decompositions in order to check the ranks of the matrices $ \wt\bB_j$ and $\wt \bC_j$
and in order to construct bi-orthogonal bases so that (\ref{t0def}) holds.
While each step of this procedure can be implemented in a numerically stable way,
the consecutive rank decisions make the aggregate procedure difficult to analyze, similar
to the case of staircase algorithms \cite{ByeMX07,DemK93a,DemK93b}. In general the strategy should be adapted
toward the goal of obtaining a realization in port-Hamiltonian form that is robust to small perturbations, see \cite{BanMNV20,MehV20} for some ways to do this.
}
\end{remark}

\subsection{Explicit solution of linear matrix inequalities via even pencils}\label{sec:Ricin}
We have seen that to check the existence of the
transformation to port-Hamiltonian form and to explicitly construct the transformation matrices $\bT,\bV$
{ in (\ref{IOtransformations})} is equivalent to considering
the solution of the linear matrix inequality (\ref{LMIcond}).
{In this subsection we combine the recursive procedure
presented in the previous subsection with the construction of a staircase like form (\cite{Van79}) for
the even pencil (\ref{evenpencil})}.

For a given real symmetric matrix $\bQ$ denote {the matrix in (\ref{LMIcond})} by
\[
\Psi_0(\bQ):=\mat{cc}-\bA^\top\bQ-\bQ\bA&\bC^\top-\bQ\bB\\
\bC-\bB^\top\bQ&\bD+\bD^\top\rix,
\]
which is supposed to be positive semidefinite.

Let $\bB_1,\bC_1$ be defined as in (\ref{feedthroughreduction}), (\ref{bcsplit}). If $\bB_1$, $\bC_1$
satisfy Part 1. of Lemma~\ref{lemgcase} then, since
$\mathsf{Ker}\,\bC_1^\top=\mathsf{Ker}\,\bB_1$, there exist real orthogonal
matrices $\wt \bU_1$, $\wt \bV_1$ (which can be obtained by performing a permuted singular value decomposition of  $\bB_1$) such that
\eq{facbc11}
\wt \bU_1^\top\bB_1\wt \bV_{1}=\mat{cc}0&0\\\boldsymbol{\Sigma}_\bB&0\rix,\quad
\wt \bV_{1}^\top\bC_1\wt \bU_1=\mat{cc}\bC_{11}&\bC_{12}\\0&0\rix
\en
where $\boldsymbol{\Sigma}_\bB$ is invertible and $\bC_{12}\boldsymbol{\Sigma}_\bB$ is real symmetric and positive definite. Transforming
the desired $\bQ$ correspondingly as
\[
\wt \bQ:=\wt \bU_1^\top\bQ\wt \bU_1=\mat{cc}\bQ_{11}& \bQ_{12}\\  \bQ_{12}^\top& \bQ_{22}\rix,
\]
then, since the linear matrix inequality (\ref{LMIcond}) implies
$\wt \bQ(\wt \bU_1^\top\bB_1)=\wt\bU_1^\top\bC_1^\top$, it follows in the transformed variables that
\eq{formx12}
 \bQ_{22}:=\bC_{12}^\top\boldsymbol{\Sigma}_\bB^{-1}>0,\quad
 \bQ_{12}:=\bC_{11}^\top\boldsymbol{\Sigma}_\bB^{-1}.
\en
It remains to determine $ \bQ_{11}$ so that
$ \wt \bQ>0$. To achieve this, we set
\begin{equation}\label{invQ22}
\bQ_0:=\mat{cc} \bQ_{12} \bQ_{22}^{-1} \bQ_{12}^\top& \bQ_{12}\\
 \bQ_{12}^\top& \bQ_{22}\rix
=\bT_0^{-\top}\mat{cc}0&0\\0&\bQ_{22}\rix\bT_0^{-1},\quad \bT_0=\mat{cc}\bI&0\\
-\bQ_{22}^{-1}\bQ_{12}^\top&\bI\rix
\end{equation}
and we clearly have that $\bQ_0\geq 0$. Then we can rewrite $\wt \bQ$ as
\[
\wt \bQ=\mat{cc} \bQ_1&0\\0&0\rix+\bQ_0
=\bT_0^{-\top}\left(\mat{cc}\bQ_1&0\\0&0\rix+\mat{cc}0&0\\0&\bQ_{22}\rix\right)\bT_0^{-1}
=\bT_0^{-\top}\mat{cc}\bQ_1&0\\0&\bQ_{22}\rix\bT_0^{-1},
\]
partitioned analogously and we obtain that $\wt \bQ>0$ if and only if $ \bQ_1>0$. Let $\bV_0$ be the orthogonal matrix given in (\ref{feedthroughreduction}).
By performing a congruence transformation on $\Psi_0(\bQ)$ with
\[
 \bZ_0=\mat{cc} \bT_0&0\\0&\bV_0\rix,\quad
\bT_0:=\wt \bU_1\bT_0,\quad
 \bV_0:= \bV_0\mat{cc}\wt \bV_{1}&0\\0&\bI\rix
\]
and using the fact that $\wt\bQ(\wt \bU_1^\top\bB_1\wt \bV_{1})=(\wt \bV_{1}^\top\bC_1\wt \bU_1)^\top$
for any real symmetric $\bQ_1$, it follows that
\[
\bZ_0^{\top}\Psi_0(\bQ)\bZ_0=\mat{ccc}
-(\bT_0^{-1}\bA\bT_0)^\top\bT_0^{\top}\bQ\bT_0
-\bT_0^{\top}\bQ\bT_0
(\bT_0^{-1}\bA\bT_0)
&0&\bT_0^{\top}\bC_2^\top-\bT_0^{\top}\bQ\bT_0
(\bT_0^{-1}\bB_2)\\
0&0&0\\
\bC_2\bT_0-(\bT_0^{-1}\bB_2)^\top
\bT_0^{\top}\bQ\bT_0&0&\bS_2\rix.
\]
Partitioning
\[
\bT_0^{-1}\bA\bT_0=:\mat{cc}\bA_{11}&\bA_{12}\\\bA_{21}&\bA_{22}\rix,\quad
\bT_0^{-1}\bB_2=\mat{c}\bB_{13}\\\bB_{23}\rix,\quad
\bC_2\bT_0=\mat{cc}\bC_{31}&\bC_{32}\rix,
\]
and using the fact
\[
\bT_0^\top\bQ\bT_0=\mat{cc}\bQ_1&0\\0&\bQ_{22}\rix,
\]
we obtain that
\[
\bZ_0^\top\Psi_0 (\bQ)\bZ_0=\mat{cccc}
-\bA_{11}^\top\bQ_1-\bQ_1\bA_{11}&-\bQ_1\bA_{12}-\bA_{21}^\top\bQ_{22}
&0&\bC_{31}^\top-\bQ_1\bB_{13}\\
-\bA_{12}^\top\bQ_1-\bQ_{22}\bA_{21}&-\bA_{22}^\top\bQ_{22}-\bQ_{22}\bA_{22}&0&\bC_{32}^\top-\bQ_{22}\bB_{23}\\
0&0&0&0\\
\bC_{31}-\bB_{13}^\top\bQ_1&
\bC_{32}-\bB_{23}^\top\bQ_{22}&0&\bS_2\rix.
\]
In this way, we have that  (\ref{LMIcond}) holds  for some $\bQ>0$ if and only if
\begin{eqnarray*}
\Psi_1(\bQ_1)&:=&
\mat{c|cc}
-\bA_{11}^\top\bQ_1-\bQ_1\bA_{11}&-\bQ_1\bA_{12}-\bA_{21}^\top\bQ_{22}
&\bC_{31}^\top-\bQ_1\bB_{13}\\\hline
-\bA_{12}^\top\bQ_1-\bQ_{22}\bA_{21}&-\bA_{22}^\top\bQ_{22}-\bQ_{22}\bA_{22}&\bC_{32}^\top-\bQ_{22}\bB_{23}\\
\bC_{31}-\bB_{13}^\top\bQ_1&
\bC_{32}-\bB_{23}^\top\bQ_{22}&{\color{black}\bS_2}\rix\\
&=:&\mat{cc}-\bA_1^\top\bQ_1-\bQ_1\bA_1
&\bC_1^\top-\bQ_1\bB_1\\
\bC_1-\bB_1^\top\bQ_1& \bD_1+\bD_1^\top\rix\ge 0
\end{eqnarray*}
holds for some real symmetric positive definite $\bQ_1$, where
\[
\bA_1=\bA_{11},\quad
 \bC_1=\mat{c}-\bQ_{22}\bA_{21}\\\bC_{31}\rix,\quad
\bB_1=\mat{cc}\bA_{12}&\bB_{13}\rix,
\]
and
\[
\bD_1+\bD_1^\top=\mat{cc}
-\bA_{22}^\top\bQ_{22}-\bQ_{22}\bA_{22}&\bC_{32}^\top-\bQ_{22}\bB_{23}\\
\bC_{32}-\bB_{23}^\top\bQ_{22}&\bS_2\rix.
\]
This construction has reduced the solution of the linear matrix inequality (\ref{LMIcond}) to  the solution
of a smaller linear matrix inequality of the same form. Thus, we can again proceed in a recursive manner with the same reduction process until either the condition in Part 1. of  Lemma~\ref{lemT} no longer holds (in which case no solution exists) or
$\bD_k+\bD_k^\top$ is positive definite for some $k$.

This reduction process can be considered as the construction of a structured staircase form for
the even pencil (\ref{evenpencil}).
By applying a congruence transformation to the pencil (\ref{evenpencil})
with the matrix
\[
\bY_0=\mat{ccc}\bT_0^{-\top}&0&0\\0&\bT_0&0\\0&0&\bV_0\rix,
\]
it follows that
\[
\lambda\mat{cc|cc|ccc}0&0&\bI&0&0&0&0\\
0&0&0&\bI&0&0&0\\\hline
-\bI&0&0&0&0&0&0\\
0&-\bI&0&0&0&0&0\\\hline
0&0&0&0&0&0&0\\
0&0&0&0&0&0&0\\
0&0&0&0&0&0&0\rix
-\mat{cc|cc|ccc} 0&0&\bA_{11}&\bA_{12}&0&0&\bB_{13}\\
0&0&\bA_{21}&\bA_{22}&\boldsymbol{\Sigma}_B&0&\bB_{23}\\\hline
\bA_{11}^\top&\bA_{21}^\top&0&0&0&0&\bC_{31}^\top\\
\bA_{12}^\top&\bA_{22}^\top&0&0&\bC_{12}^\top&0&\bC_{32}^\top\\\hline
0&\boldsymbol{\Sigma}_B^\top&0&\bC_{12}&0&0&0\\
0&0&0&0&0&0&0\\
\bB_{13}^\top&\bB_{23}^\top&\bC_{31}&\bC_{32}&0&0&\bS_2\rix.
\]
%
By performing another congruence transformation with the matrix
\[
\wt\bY_0=\mat{cc|cc|ccc}\bI&0&0&0&0&0&0\\
0&\bI&0&-\bQ_{22}&0&0&0\\\hline
0&0&\bI&0&0&0&0\\
0&0&0&\bI&0&0&0\\\hline
0&0&0&0&\bI&0&0\\
0&0&0&0&0&\bI&0\\
0&0&0&0&0&0&\bI\rix,
\]
the pencil becomes
\begin{eqnarray*}
&&\lambda\mat{cc|cc|ccc}0&0&\bI&0&0&0&0\\
0&0&0&\bI&0&0&0\\\hline
-\bI&0&0&0&0&0&0\\
0&-\bI&0&0&0&0&0\\\hline
0&0&0&0&0&0&0\\
0&0&0&0&0&0&0\\
0&0&0&0&0&0&0\rix\\
&& \qquad -\mat{cc|cc|ccc} 0&0&\bA_{11}&\bA_{12}&0&0&\bB_{13}\\
0&0&\bA_{21}&\bA_{22}&\boldsymbol{\Sigma}_B&0&\bB_{23}\\\hline
\bA_{11}^\top&\bA_{21}^\top&0&-\bA_{21}^\top\bQ_{22}&0&0&\bC_{31}^\top\\
\bA_{12}^\top&\bA_{22}^\top&-\bQ_{22}\bA_{21}&-\bA_{22}^\top\bQ_{22}-\bQ_{22}\bA_{22}
&0&0&\bC_{32}^\top-\bQ_{22}\bB_{23}\\\hline
0&\boldsymbol{\Sigma}_B^\top&0&0&0&0&0\\
0&0&0&0&0&0&0\\
\bB_{13}^\top&\bB_{23}^\top&\bC_{31}&\bC_{32}-\bB_{23}^\top\bQ_{22}&0&0& \bS_2\rix.
\end{eqnarray*}
By further moving the last block row and column to the fifth position and then the 2nd block
row and column to the fifth position, i.e., by performing a congruence permutation with
\[
\mat{cc|cc|ccc}
\bI&0&0&0&0&0&0\\
0&0&0&0&\bI&0&0\\\hline
0&\bI&0&0&0&0&0\\
0&0&\bI&0&0&0&0\\\hline
0&0&0&0&0&\bI&0\\
0&0&0&0&0&0&\bI\\
0&0&0&\bI&0&0&0\rix,
\]
one obtains
\[
\lambda \mat{ccc|ccc}
0&\bI&0&0&0&0\\
-\bI&0&0&0&0&0\\
0&0&0&-\Gamma_1&0&0\\\hline
0&0&\Gamma_1^\top&0&0&0\\
0&0&0&0&0&0\\
0&0&0&0&0&0\rix -
\mat{ccc|ccc} 0&\bA_{1}&\bB_{1}&0&0&0\\
\bA_{1}^\top&0&\bC_{1}^\top&\bA_{21}^\top&0&0\\
\bB_{1}^\top&\bC_{1}&\bD_1+\bD_1^\top&\Delta_1&0&0\\\hline
0&\bA_{21}&\Delta_1^\top&0&\boldsymbol{\Sigma}_1&0\\
0&0&0&\boldsymbol{\Sigma}_1^\top&0&0\\
0&0&0&0&0&0
\rix,
\]
where
\[
\Gamma_1=\mat{c}\bI\\0\rix,\quad\Delta_1=\mat{c}\bA_{22}^\top\\\bB_{23}^\top\rix,\quad
\boldsymbol{\Sigma}_1:=\boldsymbol{\Sigma}_B
\]
and $\bA_1,\bB_1,\bC_1$, $\bD_1+\bD_1^\top$ are as defined before.
In this way, we may repeat the reduction process on the (1,1) block, which corresponds to $\Psi_1$.
In order to exploit the block structures of the pencil we use a slightly different compression technique
for $\bD_1+\bD_1^\top$. Note that we may write
\[
\bD_1+\bD_1^\top=\mat{cc}\bD_{11}&\bD_{12}\\\bD_{12}^\top&\bS_2\rix,
\]
with $\bS_2$ symmetric positive definite. Then we have
\[
\bD_1+\bD_1^\top=\mat{cc}\bI&0\\\bS_2^{-1}\bD_{12}^\top&\bI\rix^\top
\mat{cc}\bD_{11}- \bD_{12}\bS_2^{-1}\bD_{12}^\top&0\\0&\bS_2\rix
\mat{cc}\bI&0\\\bS_2^{-1}\bD_{12}^\top&\bI\rix.
\]
Let
\[
\bD_{11}- \bD_{12}\bS_2^{-1}\bD_{12}^\top =\bZ_1\mat{cc}0&0\\0&\wt \bS_2\rix \bZ_1^\top,
\]
where $\wt\bS_2$ is invertible and $\bZ_1$ is orthogonal. Then
\[
\mat{cc}\bZ_1&0\\-\bS_2^{-1}\bD_{12}^\top\bZ_1&\bI\rix^\top
(\bD_1+\bD_1^\top)
\mat{cc}\bZ_1&0\\-\bS_2^{-1}\bD_{12}^\top\bZ_1&\bI\rix
=\mat{c|cc}0&0&0\\\hline 0&\wt\bS_2&0\\0&0&\bS_2\rix
=:\mat{cc}0&0\\0&\hat \bS_2\rix.
\]
A necessary condition for the existence of a transformation to pH form is that $\hat \bS_2>0$ or equivalently $\wt\bS_2>0$. If this holds, then using the fact that
\[
\mat{cc}\bZ_1&0\\-\bS_2^{-1}\bD_{12}^\top\bZ_1&\bI\rix^\top\Gamma_1 =\mat{c}\bZ_1^\top\\0\rix,
\]
by performing a congruence transformation on the 3rd block rows and columns
with
\[\mat{cc}\bZ_1&0\\-\bS_2^{-1}\bD_{12}^\top\bZ_1&\bI\rix\]
and another congruence transformation on the fourth block row and column with $\bZ_1$
we obtain the pencil
\[
\lambda \mat{cccc|ccc}
0&\bI&0&0&0&0&0\\
-\bI&0&0&0&0&0&0\\
0&0&0&0&-\Gamma_{11}&0&0\\
0&0&0&0&-\Gamma_{21}&0&0\\\hline
0&0&\Gamma_{11}^\top&\Gamma_{21}^\top&0&0&0\\
0&0&0&0&0&0&0\\
0&0&0&0&0&0&0\rix -
\mat{cccc|ccc} 0&\bA_{1}&\bB_{11}&\bB_{12}&0&0&0\\
\bA_{1}^\top&0&\bC_{11}^\top&\bC_{21}^\top&\Delta_{11}&0&0\\
\bB_{11}^\top&\bC_{11}&0&0&\Delta_{21}&0&0\\
\bB_{12}^\top&\bC_{21}&0&\hat \bS_2&\Delta_{31}&0&0\\\hline
0&\Delta_{11}^\top&\Delta_{21}^\top&\Delta_{31}^\top&0&\boldsymbol{\Sigma}_1&0\\
0&0&0&0&\boldsymbol{\boldsymbol{\Sigma}}_1^\top&0&0\\
0&0&0&0&0&0&0\rix,
\]
where
\[
\mat{c}\Gamma_{11}\\\Gamma_{12}\rix
=\Gamma=\mat{cc}\bI&0\\\hline 0&\bI\\0&0\rix.
\]
In order to proceed, $\bB_{11},\bC_{11}^\top$ must satisfy the same conditions as
$\bB_1,\bC_1^\top$.  If these conditions hold, then we can perform a second set of
congruence transformation and transform the pencil to
\begin{eqnarray*}
&&\lambda \mat{ccc|ccc|ccc}
0&\bI_\ell&0&0&0&0&0&0&0\\
-\bI_\ell&0&0&0&0&0&0&0&0\\
0&0&0&-\Gamma_{2}&0&0&-\wt\Gamma_{11}&0&0\\\hline
0&0&\Gamma_2^\top&0&0&0&0&0&0\\
0&0&0&0&0&0&-\wt\Gamma_{21}&0&0\\
0&0&0&0&0&0&-\wt\Gamma_{31}&0&0\\\hline
0&0&\wt\Gamma_{11}^\top&0&\wt\Gamma_{21}^\top&\wt\Gamma_{131}^\top&0&0&0\\
0&0&0&0&0&0&0&0&0\\
0&0&0&0&0&0&0&0&0\\
0&0&0&0&0&0&0&0&0\rix \\
&&-
\mat{ccc|ccc|ccc} 0&\bA_{2}&\bB_{2}&0&0&0&0&0&0\\
\bA_{2}^\top&0&\bC_{2}^\top&\Theta_{2}^\top&0&0&\wt\Delta_{11}&0&0\\
\bB_{2}^\top&\bC_{2}&\bD_2+\bD_2^\top&\Delta_2&0&0&\wt\Delta_{21}&0&0\\\hline
0&\Theta_2&\Delta_2^\top&0&\boldsymbol{\Sigma}_2&0&0&0&0\\
0&0&0&\boldsymbol{\Sigma}_2^\top&0&0&\wt\Delta_{31}&0&0\\
0&0&0&0&0&0&\wt\Delta_{41}&0&0\\\hline
0&\wt\Delta_{11}^\top&\wt\Delta_{21}^\top&0&\wt\Delta_{31}^\top&\wt\Delta_{41}^\top&0&\boldsymbol{\Sigma}_1&0\\
0&0&0&0&0&0&\boldsymbol{\boldsymbol{\Sigma}}_1^\top&0&0\\
0&0&0&0&0&0&0&0&0\rix,
\end{eqnarray*}
where
\[\wt \Gamma_{11}=\mat{c}0\\\Gamma_{21}\rix,\quad
\mat{c}
\wt\Gamma_{11}\\\wt\Gamma_{21}\rix
=\Gamma_{31},\quad \Gamma_2=\mat{c}\bI\\0\rix,
\]
and with the partitioning $\Delta_{11}=\mat{c}\Delta_{11,1}\\\Delta_{11,2}\rix$,
\[
\wt \Delta_{11}=\Delta_{11,1},\quad
\wt \Delta_{21}=\mat{c}\Delta_{11,2}\\\Delta_{31}\rix, \quad
\mat{c}\wt \Delta_{31}\\\wt\Delta_{41}\rix=\Delta_{21}.
\]
This reduction process continues as long as all the required conditions hold, until for some $k$, $\bD_k+\bD_k^\top>0$.
If this is the case, then the pencil (\ref{evenpencil}) is reduced to  an even pencil that has the eigenvalue $\infty$ with equal algebraic and geometric multiplicity (it is of index one)
\begin{equation}\label{hampen}
\lambda\mat{ccc}0&\bI_\ell&0\\-\bI_\ell&0&0\\0&0&0\rix
-\mat{ccc}0&\bA_k&\bB_k\\
\bA_k^\top&0&\bC_k^\top\\
\bB_k^\top&\bC_k&\bD_k+\bD_k^\top\rix.
\end{equation}
Note the above process is actually a special staircase form reduction process that deflates
the singular part and higher index of the eigenvalue infinity of the even pencil (\ref{evenpencil}), \cite{BruM08,ByeMX07}.

{In order to determine a solution $\bQ=\bT^\top\bT>0$, we still need to use the method developed in the
previous section to find $\bQ_k=\bT_k^\top\bT_k>0$ for solving the
KYP inequality (\ref{LMIcond}) with $(\bA_k,\bB_k,\bC_k,\bD_k)$.

If we consider the positive definite solution of the KYP equation,
we can directly check if (\ref{hampen}) has a deflating subspace associated with a set of $\ell$ finite eigenvalues chosen such that the deflating subspace is as in \eqref{evenpencil}. If such a deflating subspace
exists and $\bW_1$ is invertible, then
we can compute a Hermitian positive definite matrix $\bQ_k$ associated with \eqref{hampen}. We comment that an analogous characterization of solutions of singular H$^\infty$ control problems via matrix pencils was given in \cite{CopS92}.
}

\begin{remark}\label{rem:notepass}{\rm
Note that to check the passivity of (\ref{GenSys}) it is only necessary to have a positive
semidefinite solution to (\ref{LMIcond}). Thus, if one only wants to check passivity, then {Part a) in Lemma~\ref{lemT}}
can be relaxed to $\mathsf{Ker}\,\bB_1\subseteq \mathsf{Ker}\,\bC_1^\top$. In this case the transformation to the form
(\ref{facbc11}) can still be made, but   $ \bQ_{22}$ in (\ref{formx12}) is only
positive semidefinite, and $\mathsf{Ker}\,\bQ_{22}\subseteq\mathsf{Ker}\,\bQ_{12}$.
Then $\bQ_0$ can still be defined but instead of $ \bQ_{11}^{-1}$ one needs to use
the Moore-Penrose pseudoinverse, see \cite{GolV96}, of $ \bQ_{11}$. However,
in this case $\bQ_0$ and the resulting solution $\bQ$ cannot be
positive definite.
Thus, in this situation, (\ref{GenSys}) may be a passive
system that cannot be transformed to a standard port-Hamiltonian system of the form \eqref{FullpHSys} system.

A simple example
is the scalar system
\[
\dot\bz =-\bz+2\bu,\quad  \by=0\bz+0\bu.
\]
This system is passive (but not strictly passive) and the matrix inequality \eqref{LMIcond} has the unique singular solution $\bQ=0$.
So this system cannot be transformed to a port-Hamiltonian system of the form~\eqref{FullpHSys}.
%
%
%
%
In this case then one has to use a descriptor formulation~\eqref{FulldaepHSys}.
}
\end{remark}

To illustrate the analysis procedures, consider the following example.
\begin{example}
{\rm In the finite element analysis of disc brake squeal \cite{GraMQSW16}, the model is a very large-scale
second-order system of differential equations with approximately a million degrees of freedom, that furthermore also depends on parameters, e.g., the disc speed $\omega$. If no further constraints are incorporated, then in the stationary case the system takes the form
\[
\bM\ddot{{\bq}}+\left(\bC_1+\frac{\omega_{r}}{\omega}\bC_R+\frac{\omega}{\omega_{r}} \bC_G\right)\dot{{\bq}}+\left(\bK_1+\bK_R+\left(\frac{\omega}{\omega_{r}}\right)^2 \bK_{G}\right){\bq}=\bB \bu, \quad\by=\bB^\top \bq\ ,
\]
where $\bM=\bM^\top>0$ is the mass matrix,
$\bC_{1}=\bC_1^\top\geq 0$ models material damping, $\bC_{G}=-\bC_G^\top$ models  gyroscopic effects,
$\bC_{R}=\bC_R^\top\geq 0$ models friction induced damping, $\bK_{1}=\bK_1^\top>0$ is the stiffness matrix,
$\bK_{R}=\bK_2+\bN$ with $\bK_2=\bK_2^\top$ and $\bN=-\bN^\top$, is a nonsymmetric matrix modeling circulatory effects, $\bK_{G}=\bK_G^\top\geq 0$ is the geometric stiffness matrix, and $\omega$ is the rotational speed of the disc with reference velocity   $\omega_r$. In industrial brake models, the matrices $\bD:=\bC_1+\frac{\omega_{r}}{\omega}\bC_R$, and $\bN$ are sparse and have very low rank  (approx. $2000$) corresponding to finite element nodes associated with the brake pad. Setting $\bG:=\frac{\omega}{\omega_{r}}\bC_G$, $\bK=\bK_1+\bK_2+(\frac{\omega}{\omega_{r}})^2\bK_G$, we may assume that $\bK>0$. Here in the practical design  a \emph{shim} is attached to the brake pad which may be interpreted as choosing the input as output feedback $\bu=\bD_\bS \bB^\top\bq$ in order to stabilize the system in a given range of disk speeds.

Then, introducing $\bp=\bM \dot \bq$, we can write the system in first order form
\[
\mat{c} \dot \bp \\ \dot \bq \rix = (\bJ-\bR)\bQ \mat{c} \bp \\ \bq \rix,
\]
where
\begin{eqnarray*}
\bJ&:=&\mat{cc} -\bG  &  -(\bI+\frac 12(\bN-\bB\bD_{\bS}\bB^\top) \bK^{-1})\\
(\bI+\frac 12(\bN-\bB\bD_\bS\bB^\top)\bK^{-1})^\top & 0 \rix,\\
\bR&:=&\mat{cc} \bD & \frac 12(\bN-\bB\bD_{\bS}\bB^\top) \bK^{-1}\\
(\frac 12(\bN-\bB\bD_{\bS}\bB^\top) \bK^{-1})^\top& 0 \rix,\quad
\bQ=\mat{cc}\bM^{-1}&0\\0&\bK\rix.
\end{eqnarray*}
Regardless of the choice of {the matrix $\bD_\bS$ in the feedback $\bu$}, the matrix $\bR$
is indefinite as long as $\bN\neq 0$ (then $\bN-\bB\bD_\bS\bB^\top\neq 0$); thus, it is clear that for this system we cannot read off its stability and it is definitely unstable if
$\bx^H \bR \bx<0$ for some eigenvector $\bx$ of $\bJ$. {This shows that the rewriting of $\bA$ as $(\bJ -\bR) \bQ)$ alone will not be enough to check stability, but if the system is asymptotically stable then a further transformation as described here will be necessary.}
}
\end{example}
\begin{remark}\label{rem:dae1}{\rm
Many of the results described in this work can be extended to the case of general descriptor systems having the form~\eqref{pHdaeelem} or~\eqref{FulldaepHSys} but with singular $\bE$. This is an open problem and currently a topic of active research.}
\end{remark}

\begin{remark}\label{rem:extreme}
{\rm In this paper, we only consider positive definite solutions $\bQ$ of the KYP matrix inequality~\eqref{LMIcond}. In principle, this condition could be relaxed to requiring that  $\bQ$ is only positive semidefinite. However, several further subtle problems arise in this case, and many of them lead to the need to consider descriptor systems.
See \cite{CheGH22,MehMW18,MehS23} for a detailed  discussion.}
\end{remark}
\begin{remark}\label{rem:robust}{\rm
In view of the fact that the KYP matrix inequality \eqref{LMIcond}
may have many solutions, one may choose the solution so that the resulting pH representation is robust to data or numerical errors. To do this one could choose a robustness measure like to the distance to instability or non-passivity or to consider a representation that is far away from the boundary of the solution set of \eqref{LMIcond}.
Partial results in this direction have been obtained in \cite{BanMNV20,BeaMV19,MehV20a}.
}
\end{remark}

\section{Numerical Considerations}\label{sec:numerical}
In this section we discuss some numerical issues that arise when implementing the procedure as discussed in the last section.  An associated \textsc{matlab} script is available for download on the \textsc{MathWorks} FileExchange (under \verb!nearby-pH-realization!)  and at GitHub\footnote{\texttt{https://github.com/christopherbeattie/nearby-pH-realization}}.

Most steps of our procedure are implemented in a straightforward way using standard techniques from numerical linear algebra that ensure backward stability.
There are a few places where forward stability may be lost and large relative errors could occur. In particular, the similarity transformations with respect to $\bT_i$ in \eqref{kstepineq} as well as the linear solves implicit in the inversions of $Q_{22}$ in \eqref{invQ22} used in forming the matrices $\bT_i$. If these are ill-conditioned then large relative errors may arise.

Another difficulty is the non-uniqueness of solutions to the linear matrix inequality \eqref{LMIred}. To make the solution unique one can optimize a quality measure like the distance to instability or non-passivity, or try to find the analytic center of the solution, see \cite{BanMNV20,MehV20}. All these are difficult and expensive optimization problems on top of all the computational work that has to be carried out.
How to do this efficiently is an open question, even for the case that the system is minimal and $\hat \bS$ is well-conditioned with respect to inversion. The difficulty arises, in particular, since the solution sets of \eqref{LMIred} and \eqref{RIE}  are rather difficult to characterize, see \cite{MehX24} for a detailed analysis based on eigenvalue perturbation theory.

\begin{remark}\label{rem:newton}{\rm
Note that if the skew-symmetric/symmetric pencil \eqref{hampen} has purely imaginary eigenvalues, then the solution of the Riccati equation associated with \eqref{RIE} can only be computed with the Newton method of \cite{Ben97} which has been implemented e.g.  in \cite{BenW20}. In this case none of the usual approaches utilizing invariant subspaces will be fully satisfactory and even the Newton method might not achieve a quadratic rate of convergence, displaying only linear convergence \cite{GuoL98}.
}
\end{remark}

We have tested our procedure for a large number of examples with randomly generated stable and passive systems which were produced from  pH systems by multiplying out the factors. In each example the procedure yielded the same  pH representation. 

\begin{remark}\label{rem:modelphilosopy}{\rm
As we have discussed in Section~\ref{sec:suff}, we reduce the problem to a subproblem, where $\hat \bC$ is invertible and thus where $\eqref{RIE}$ can be formed and is solvable. For this we have to make several rank decisions or regularization steps which result in small perturbations. This is often justified, since the coefficient matrices $\bA$, $\bB$, $\bC$, $\bD$ are typically not exact because they arise often from a data based realization, interpolation or model reduction process. So we can make small perturbations to these data to regularize the problem if this does not  change the resulting (hopefully) robust pH representation.

This is common practice when solving \eqref{LMIred}, see e.g. \cite{AndM74,HilM80}, where often $\hat \bS$ is perturbed to be invertible, so that the pencil \eqref{evenpencil} is regular and of index at most one, see Remark~\ref{rem_even}. However, as discussed in Section~\ref{sec:suff}, see Remark~\ref{rem:almostsingular}, since we are able to deal with a singular $\hat \bS$ we perform the regularization in a different way. For this we have to perform the rank decisions in \eqref{feedthroughreduction}, \eqref{factbc1}, \eqref{facbc11}, which are critical in the process.
As in most staircase algorithms \cite{Van79}, it is recommended to make conservative decisions, i.e., to assume smaller rank if the decision is difficult using the usual rank decision procedures \cite{DemK93a,DemK93b}.
}
\end{remark}

\begin{remark}\label{rem:nonunique2}{\rm Note that if there is no solution to the  linear matrix inequality \eqref{LMIcond} our procedure will produce a positive definite solution $Q$ to a slightly perturbed linear matrix inequality for which a solution can be assured, see Example~\ref{ex:smallex} below.}
\end{remark}

\begin{remark}\label{rem:nonunique}{\rm Since the solution to \eqref{LMIcond} is not unique, even if we start with a pH realization of a system, i.e., when the system is in the form \eqref{pHdef} with quadratic Hamiltonian $\mathcal H=\frac 12 x^\top Qx$, our procedure in general, may not return the same representation, see e.g. Example~\ref{ex:smallex} below.}
\end{remark}

\begin{example}\label{ex:smallex}{\rm
Consider a system \eqref{GenSys} with $\bA=(\bJ-\bR) \bQ$, $\bB= \bF-\bP$,\ $\bC=(\bF+\bP)^\top$,\ $\bD=\bS+\bN$, where
$\bJ=0$, $\bN=0$, $\bQ=\bI_4$,
\[\bR=\begin{bmatrix} 1 & &&\\ & 2 && \\ && 3 & \\ &&& 4 \end{bmatrix},\ \bB=\begin{bmatrix} 1 & 2\\ 1 & -1  \\ 1& 3 \\ 1 & 1 \end{bmatrix},\
\bP=\frac 14 \begin{bmatrix}\epsilon & -1\\ \epsilon & 1  \\ \epsilon & -1 \\ \epsilon & 1 \end{bmatrix},\ \bS=\frac 12\begin{bmatrix} \epsilon & 0\\ 0 & 1 \end{bmatrix}.
\]

We ran our script with different values of $\epsilon=0, \pm 2.*e{-51},\ldots,\pm 2.*e{-9}$. When $\epsilon$ is very close to $0$ but negative, the system is slightly non-passive, since  the matrix $\bW$ has a negative eigenvalue. Our script always returns a positive definite $\bQ$ that changes only slightly as  $\epsilon$ varies (regardless of sign), while the inequality \eqref{LMIred} remains semidefinite but numerically singular with four very small nonzero eigenvalues. For  $\epsilon=0,\pm 10^{-16}$  the pencil \eqref{evenpencil} has index $\geq 2$, otherwise the index is $1$.  Nonetheless, our procedure  works in all cases, see Figure~\ref{figures_ex_nonunique}, {where by 'residual' we denote the norm of the perturbation that we have to add to the linear matrix inequality \eqref{LMIcond} to guarantee the existence of a port-Hamiltonian formulation}.
\begin{figure}
		\includegraphics[width = 0.5\textwidth]{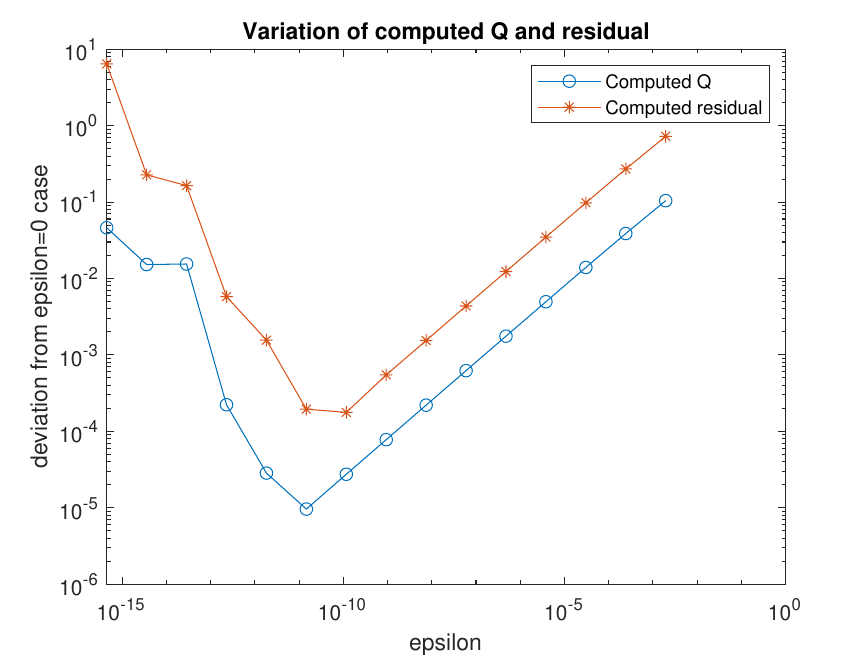}\qquad
	\includegraphics[width = 0.5\textwidth]{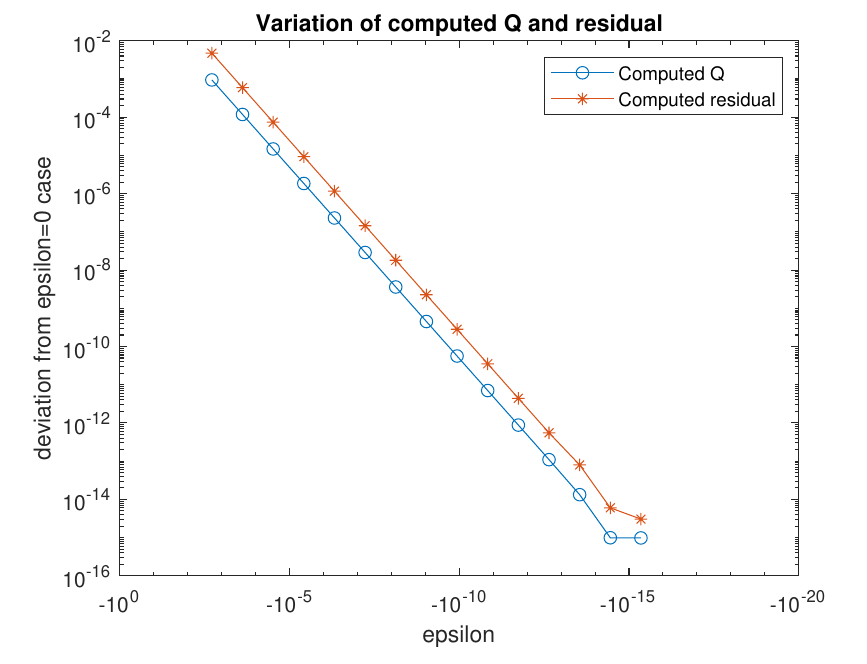}
	\caption{Deviation of computed $\bQ$ and residual of \eqref{LMIcond} from solution for $\epsilon=0$ for varying $\epsilon$.}
	\label{figures_ex_nonunique}
\end{figure}
Running our script for a system \eqref{GenSys} with the same coefficients and $\epsilon$ values as in Example~\ref{ex:smallex},
but  $\bQ$ random positive definite, we obtain the results in
Figure~\ref{figures_ex_randomQ}.
\begin{figure}
		\includegraphics[width = 0.5\textwidth]{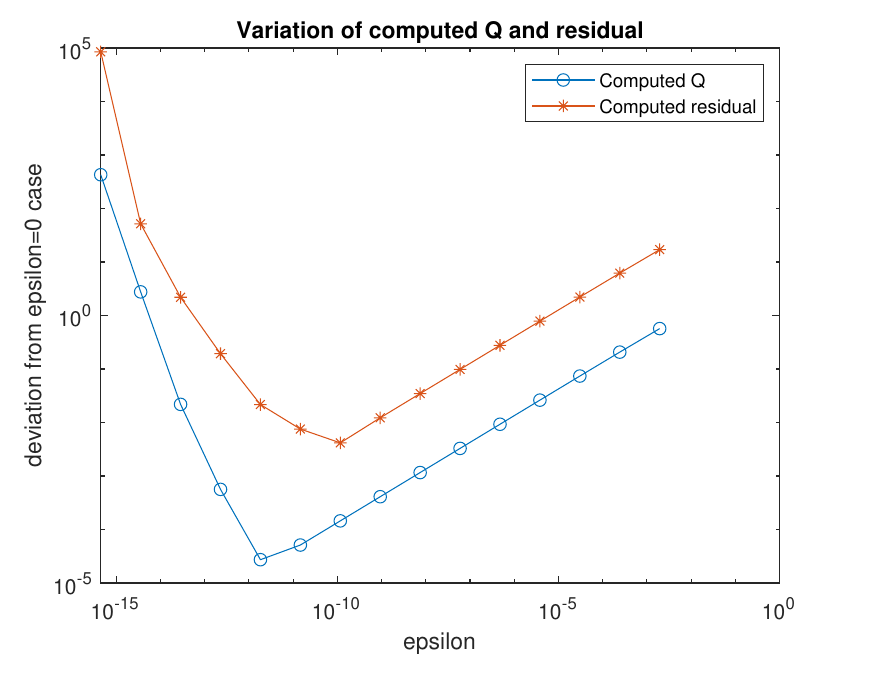}\qquad
	\includegraphics[width = 0.5\textwidth]{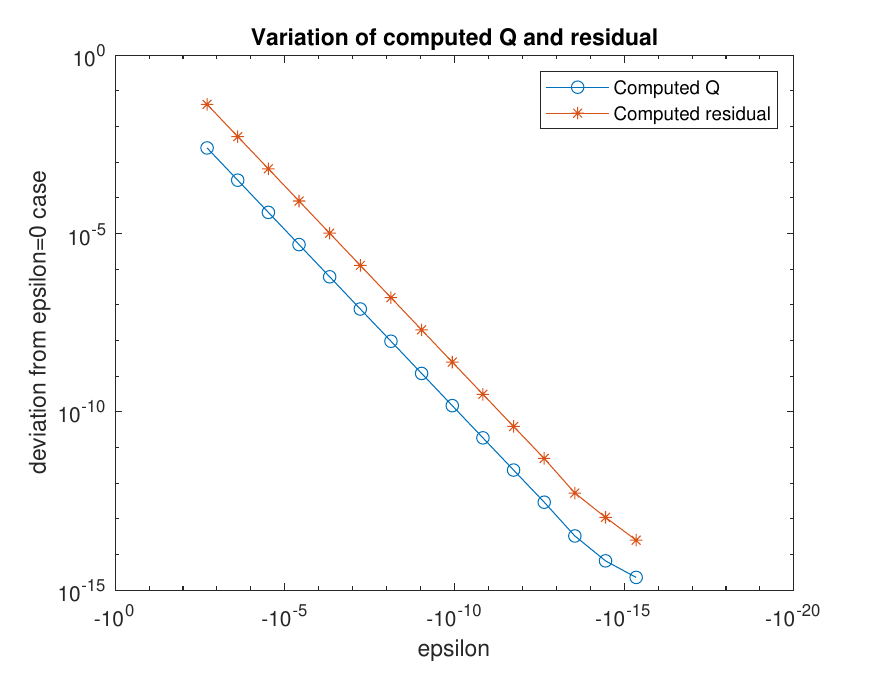}
	\caption{Deviation of $\bQ$ and residual for \eqref{LMIcond} from solution of $\epsilon=0$ for varying $\epsilon$ and random positive definite $\bQ$.}
	\label{figures_ex_randomQ}
\end{figure}
}
\end{example}

\begin{remark}\label{rem:staircase}
    {\rm
In Section~\ref{sec:suff} we have developed a method that allows to identify the solvability of the matrix inequality
\eqref{LMIcond} via the computation of a staircase like form and deflating subspaces of the even matrix pencil \eqref{evenpencil}.
In numerical practice, to avoid the use of consecutive rank decisions, one may apply a so called  derivative array approach, see \cite{KunM24}, where the parts associated with the finite eigenvalues, the part associated with the infinite eigenvalues and the singular part are separated by one sided transformations of \eqref{evenpencil} from the left generated based on an extended pencil associated with a derivative array of the DAE associated with the even pencil.
For even pencils this follows directly from a procedure developed in \cite{KunMS13} for time varying DAEs and gives a reduced system that is associated with an even pencil and an unstructured part associated with the eigenvalue $\infty$. We do not present this approach here, see \cite{KunMS13,KunM23} for details.
}
\end{remark}

\section{Conclusions}
{ Building on several well-known  results, see Section~\ref{sec:lyaric}, we have in this paper extended the characterization of
when a system is equivalent to a port-Hamiltonian system to the case of general non-minimal systems and to the case that the symmetric part of the feedthrough matrix is singular.
We have presented an explicit  procedure for the construction of the transformation matrices, and have provided an implementation. The method presented works in all tested and synthetically constructed problems.
By generating minimal perturbations to system coefficients, our procedure can also be used on a system that is not stable or not passive to produce a pH representation of a nearby system that is both stable and passive.}

Open problems include the question how to parameterize the positive definite solutions of \eqref{LMIcond} in terms of eigenvalues or pseudospectra of
the system matrix 
and the choice of adequate robustness measures to select an optimal solution of \eqref{LMIcond}.
The extension of the approach to differential-algebraic equations is another important research topic.

\subsection*{Declaration of competing interest}

The authors declare that they have no known competing financial interests or personal
relationships that could have appeared to influence the work reported in this paper.


\begin{thebibliography}{10}

\bibitem{AlaBKMM11}
R.~Alam, S.~Bora, M.~Karow, V.~Mehrmann, and J.~Moro.
\newblock Perturbation theory for {H}amiltonian matrices and the distance to
  bounded-realness.
\newblock {\em {SIAM} J. Matrix Anal. Appl.}, 32:484--514, 2011.

\bibitem{AlpL11}
D.~Alpay and I.~Lewkowicz.
\newblock The positive real lemma and construction of all realizations of
  generalized positive rational functions.
\newblock {\em Systems \& Control Let.}, 60(12):985--993, 2011.

\bibitem{AndM74}
B.~Anderson and P.~Moylan.
\newblock Synthesis of linear time-varying passive networks.
\newblock {\em IEEE Trans.  Circuits  Systems}, 21(5):678--687,
  1974.

\bibitem{BanMNV20}
D.~Bankmann, V.~Mehrmann, Y.~Nesterov, and P.~Van Dooren.
\newblock Computation of the analytic center of the solution set of the linear
  matrix inequality arising in continuous- and discrete-time passivity
  analysis.
\newblock {\em Vietnam J. Mathematics}, 48:633--660, 2020.
\newblock {\tt http://arxiv.org/abs/1904.08202}.

\bibitem{BeaG11}
C.~Beattie and S.~Gugercin.
\newblock Structure-preserving model reduction for nonlinear port-{H}amiltonian
  systems.
\newblock In {\em 50th IEEE Conf. .Decision and Control and European
  Control Conf. (CDC-ECC), 2011}, pages 6564--6569. IEEE, 2011.

\bibitem{BeaMV19}
C.~Beattie, V.~Mehrmann, and P.~{Van Dooren}.
\newblock Robust port-{H}amiltonian representations of passive systems.
\newblock {\em Automatica}, 100:182--186, 2019.

\bibitem{BeaMXZ18}
C.~Beattie, V.~Mehrmann, H.~Xu, and H.~Zwart.
\newblock Port-{H}amiltonian descriptor systems.
\newblock {\em Math. Control Signals Syst.}, 30:1--27, 2018.
\newblock https://doi.org/10.1007/s00498-018-0223-3.

\bibitem{Ben97}
P.~Benner.
\newblock Numerical solution of special algebraic {R}iccati equations via an
  exact line search method.
\newblock In {\em 1997 European Control Conf. (ECC)}, pages 3136--3141.
  IEEE, 1997.

\bibitem{BenBMX99}
P.~Benner, R.~Byers, V.~Mehrmann, and H.~Xu.
\newblock Numerical methods for linear-quadratic and {$\mathcal{H}_\infty$}
  control problems.
\newblock In G.~Picci and D.S. Gilliam, editors, {\em Dynamical Systems,
  Control, Coding, Computer Vision: New Trends, Interfaces, and Interplay},
  volume~25 of {\em Progress in Systems and Control Theory}, pages 203--222.
  Birkh{\"a}user, Basel, 1999.

\bibitem{BenLMV15}
P.~Benner, P.~Losse, V.~Mehrmann, and M.~Voigt.
\newblock Numerical linear algebra methods for linear differential-algebraic
  equations: A survey.
\newblock In {\em DAE Forum III}, pages 117--165. Springer, 2015.

\bibitem{BenW20}
P. Benner and S.W.R. Werner.
\newblock {MORLAB}--a model order reduction framework in {MATLAB} and {O}ctave.
\newblock In {\em Intern. Cong.on Mathematical Software}, pages
  432--441. Springer, 2020.

\bibitem{BoyEFB94}
S.~Boyd, L.~El Ghaoui, E.~Feron, and V.~Balakrishnan.
\newblock {\em Linear Matrix Inequalities in Systems and Control Theory}.
\newblock SIAM, Philadelphia, 1994.

\bibitem{BruM08}
T.~Br{\"u}ll and V.~Mehrmann.
\newblock {STCSSP}: A {FORTRAN} 77 routine to compute a structured staircase
  form for a (skew-)symmetric/(skew-)symmetric matrix pencil.
\newblock Preprint 31-2007, {I}nstitut f{\"u}r {M}athematik, {TU} {B}erlin,
  2007.

\bibitem{BruS13}
T.~Br\"ull and C.~Schr\"oder.
\newblock Dissipativity enforcement via perturbation of para-{H}ermitian
  pencils.
\newblock {\em IEEE Trans. Circuits Syst. I. Regul. Pap.}, 60(1):164--177,
  2013.

\bibitem{ByeMX07}
R.~Byers, V.~Mehrmann, and H.~Xu.
\newblock A structured staircase algorithm for skew-symmetric/symmetric
  pencils.
\newblock {\em Electron. Trans. Numer. Anal.}, 26:1--13, 2007.

\bibitem{ByrI03}
C.I. Byrnes and A.~Isidori.
\newblock Limit sets, zero dynamics, and internal models in the problem of
  nonlinear output regulation.
\newblock {\em IEEE Trans. Automatic Control}, 48(10):1712--1723,
  2003.

\bibitem{ByrIW91}
C.I. Byrnes, A.~Isidori, and J.C. Willems.
\newblock Passivity, feedback equivalence, and the global stabilization of
  minimum phase nonlinear systems.
\newblock {\em IEEE Trans. Autom. Control}, 36:1228--1240, 1991.

\bibitem{CheGH22}
K.~Cherifi, H.~Gernandt, and D.~Hinsen.
\newblock The difference between port-{H}amiltonian, passive and positive real
  descriptor systems.
\newblock {\em Math. Control Signals Syst.}, pages 1--23, 2023.
\newblock https://arxiv.org/pdf/2204.04990.pdf.

\bibitem{CopS92}
B.R. Copeland and M.G. Safonov.
\newblock A generalized eigenproblem solution for singular {$\mathcal{H}_2$}
  and {$\mathcal{H}_\infty$} problems.
\newblock In {\em Robust control system techniques and applications, Part 1},
  volume~50 of {\em Control Dynam. Systems Adv. Theory Appl.}, pages 331--394.
  Academic Press, San Diego, CA, 1992.

\bibitem{DemK93a}
J.~Demmel and B.~K\aa{}gstr\"om.
\newblock The generalized {S}chur decomposition of an arbitrary pencil
  {$A-zB$}: {R}obust software with error bounds and applications. {P}art {I}:
  {T}heory and algorithms.
\newblock {\em ACM Trans. Math. Softw.}, 19(2):160--174, 1993.

\bibitem{DemK93b}
J.~Demmel and B.~K\aa{}gstr\"om.
\newblock The generalized {S}chur decomposition of an arbitrary pencil
  {$A-zB$}: {R}obust software with error bounds and applications. {P}art {II}:
  {S}oftware and applications.
\newblock {\em ACM Trans. Math. Softw.}, 19(2):175--201, 1993.

\bibitem{FauMPSW22}
T.~Faulwasser, B.~Maschke, F.~Philipp, M.~Schaller, and K.~Worthmann.
\newblock Optimal control of port-hamiltonian descriptor systems with minimal
  energy supply.
\newblock {\em {SIAM} J. Cont. Optim.}, 60(4):2132--2158, 2022.

\bibitem{FreMX02}
G.~Freiling, V.~Mehrmann, and H.~Xu.
\newblock Existence, uniqueness and parametrization of {L}agrangian invariant
  subspaces.
\newblock {\em {SIAM} J. Matrix Anal. Appl.}, 23:1045--1069, 2002.

\bibitem{GolSBM03}
G.~Golo, A.J.~{van der} Schaft, P.C. Breedveld, and B.M. Maschke.
\newblock {H}amiltonian formulation of bond graphs.
\newblock In A.~Rantzer R.~Johansson, editor, {\em Nonlinear and Hybrid Systems
  in Automotive Control}, pages 351--372. Springer, Heidelberg, 2003.

\bibitem{GolV96}
G.~H. Golub and C.~F. {Van~Loan}.
\newblock {\em Matrix Computations}.
\newblock Johns Hopkins Univ. Press, Baltimore, 3rd edition, 1996.

\bibitem{GraMQSW16}
N.~Gr{\"a}bner, V.~Mehrmann, S.~Quraishi, C.~Schr\"oder, and U.~{von W}agner.
\newblock Numerical methods for parametric model reduction in the simulation of
  disc brake squeal.
\newblock {\em Z. Angew. Math. Mech.}, 96(DOI:
  10.1002/zamm.201500217):1388--1405, 2016.

\bibitem{Gri04}
S.~Grivet-Talocia.
\newblock Passivity enforcement via perturbation of {H}amiltonian matrices.
\newblock {\em IEEE Trans. Circuits and Systems}, 51:1755--1769, 2004.

\bibitem{GugPBS12}
S.~Gugercin, R.V. Polyuga, C.~Beattie, and A.J.~{van der} Schaft.
\newblock Structure-preserving tangential interpolation for model reduction of
  port-{H}amiltonian systems.
\newblock {\em Automatica}, 48:1963--1974, 2012.

\bibitem{GuoL98}
C.-H. Guo and P.~Lancaster.
\newblock Analysis and modificaton of {N}ewton’s method for algebraic
  {R}iccati equations.
\newblock {\em Math. Comp.}, 67(223):1089--1105, 1998.

\bibitem{HilM80}
D.~J. Hill and P.~J. Moylan.
\newblock Dissipative {{Dynamical Systems}}: {{Basic Input-Output}} and {{State
  Properties}}.
\newblock {\em Journal  Franklin Inst.}, 309(5):327--357, 1980.

\bibitem{HinP05}
D.~Hinrichsen and A.~J. Pritchard.
\newblock {\em Mathematical Systems Theory {I}. Modelling, State Space
  Analysis, Stability and Robustness}.
\newblock Springer-Verlag, New York, NY, 2005.

\bibitem{HorJ85}
R.~A. Horn and C.~R. Johnson.
\newblock {\em Matrix Analysis}.
\newblock Cambridge University Press, Cambridge, 1985.

\bibitem{JacZ12}
B.~Jacob and Zwart H.
\newblock {\em Linear port-{H}amiltonian systems on infinite-dimensional
  spaces}.
\newblock Operator Theory: Advances and Applications, 223.
  Birkh{\"a}user/Springer Basel AG, Basel CH, 2012.

\bibitem{Kai80}
T.~Kailath.
\newblock {\em Linear systems}, volume~1.
\newblock Prentice-Hall Englewood Cliffs, NJ, 1980.

\bibitem{Kle13}
C.~Kleijn.
\newblock {\em 20-sim 4C 2.1 Reference manual}.
\newblock Controlab Products B.V., 2013.

\bibitem{KotA10}
N.~Kottensette and P.J. Antsaklis.
\newblock Relationships between positive real, passive, dissipative and
  positive systems.
\newblock In {\em American Control Conf.}, Baltimore, MD, USA, 2010.

\bibitem{KunM23}
P.~Kunkel and V.~Mehrmann.
\newblock Local and global canonical forms for differential-algebraic equations
  with symmetries.
\newblock {\em Vietnam J. Mathematics}, 51:177--198, 2023.

\bibitem{KunM24}
P.~Kunkel and V.~Mehrmann.
\newblock {\em Differential-Algebraic Equations. Analysis and Numerical
  Solution}.
\newblock EMS Press, Berlin, Germany, 2nd edition, 2024.

\bibitem{KunMS13}
P.~Kunkel, V.~Mehrmann, and L.~Scholz.
\newblock Self-adjoint differential-algebraic equations.
\newblock {\em Math. Control Signals Syst.}, 26:47--76, 2014.

\bibitem{LanR95}
P.~Lancaster and L.~Rodman.
\newblock {\em The Algebraic {R}iccati Equation}.
\newblock Oxford University Press, Oxford, 1995.

\bibitem{LanT85}
P.~Lancaster and M.~Tismenetsky.
\newblock {\em The Theory of Matrices}.
\newblock Academic Press, Orlando, 2nd edition, 1985.

\bibitem{MasSB92}
B.M. Maschke, A.J.~{van der} Schaft, and P.C. Breedveld.
\newblock An intrinsic {H}amiltonian formulation of network dynamics:
  non-standard poisson structures and gyrators.
\newblock {\em J. Franklin Inst.}, 329:923ñ--966, 1992.

\bibitem{MehMW18}
C.~{Mehl}, V.~{Mehrmann}, and M.~{Wojtylak}.
\newblock Linear algebra properties of dissipative {H}amiltonian descriptor
  systems.
\newblock {\em {SIAM} J. Matrix Anal. Appl.}, 39:1489--1519, 2018.

\bibitem{Meh91}
V.~Mehrmann.
\newblock {\em The Autonomous Linear Quadratic Control Problem, Theory and
  Numerical Solution}, volume 163 of {\em Lecture Notes in Control and Inform.
  Sci.}
\newblock Springer-Verlag, Heidelberg, July 1991.

\bibitem{MehV20a}
V.~Mehrmann and P.~Van Dooren.
\newblock Optimal robustness of discrete-time passive systems.
\newblock {\em IMA J. of Math. Control and Inf.}, 37:1248--1269, 2020.

\bibitem{MehV20}
V.~Mehrmann and P.~Van Dooren.
\newblock Optimal robustness of port-{H}amiltonian systems.
\newblock {\em {SIAM} J. Matrix Anal. Appl.}, 41:134--151, 2020.

\bibitem{MehM19}
V.~Mehrmann and R.~Morandin.
\newblock Structure-preserving discretization for port-{H}amiltonian descriptor
  systems.
\newblock In {\em 58th IEEE Conf. Decision and Control (CDC),
  9.-12.12.19, Nice}, pages 6863--6868. IEEE, 2019.

\bibitem{MehU23}
V.~Mehrmann and B.~Unger.
\newblock Control of port-{H}amiltonian differential-algebraic systems and
  applications.
\newblock {\em Acta Numerica}, pages 395--515, 2023.

\bibitem{MehS23}
V.~Mehrmann and A.J. van~der Schaft.
\newblock Differential-algebraic systems with dissipative {H}amiltonian
  structure.
\newblock {\em Math. Control Signals Syst.}, pages 1--44, 2023.

\bibitem{MehX00}
V.~Mehrmann and H.~Xu.
\newblock Numerical methods in control.
\newblock {\em J. Comput. Appl. Math.}, 123:371--394, 2000.

\bibitem{MehX24}
V.~{Mehrmann} and H.~{Xu}.
\newblock Eigenstructure perturbations for a class of {H}amiltonian matrices
  and solutions of related {R}iccati inequalities.
\newblock {\em {SIAM} J. Matrix Anal. Appl.}, 45:1335--1360, 2024.

\bibitem{OrtSMM01}
R.~Ortega, A.J.~{van der} Schaft, Y.~Mareels, and B.M. Maschke.
\newblock Putting energy back in control.
\newblock {\em Control Syst. Mag.}, 21:18--ñ33, 2001.

\bibitem{OrtSME02}
R.~Ortega, A.J.~{van der} Schaft, B.M. Maschke, and G.~Escobar.
\newblock Interconnection and damping assignment passivity-based control of
  port-controlled {H}amiltonian systems.
\newblock {\em Automatica}, 38:585--596, 2002.

\bibitem{PaiV81}
C.C. Paige and C.F. {Van~Loan}.
\newblock A {S}chur decomposition for {H}amiltonian matrices.
\newblock {\em Linear Algebra Appl.}, 14:11--32, 1981.

\bibitem{PolS10}
R.~V. Polyuga and A.J.~{van der} Schaft.
\newblock Structure preserving model reduction of port-{H}amiltonian systems by
  moment matching at infinity.
\newblock {\em Automatica}, 46:665--672, 2010.

\bibitem{Sch04}
A.J.~{van der} Schaft.
\newblock Port-{H}amiltonian systems: network modeling and control of nonlinear
  physical systems.
\newblock In {\em Advanced Dynamics and Control of Structures and Machines},
  {CISM} Courses and Lectures, Vol. 444. Springer Verlag, New York, N.Y., 2004.

\bibitem{Sch06}
A.J.~{van der} Schaft.
\newblock Port-{H}amiltonian systems: an introductory survey.
\newblock In J.L.~Verona M.~Sanz-Sole and J.~Verdura, editors, {\em Proc. of
  the Intern. Congr. Mathematicians, vol. III, Invited Lectures},
  pages 1339--1365, Madrid, Spain, 2006.

\bibitem{SchM95}
A.J.~{van der} Schaft and B.M. Maschke.
\newblock The {H}amiltonian formulation of energy conserving physical systems
  with external ports.
\newblock {\em Arch. Elektron. {\"U}bertragungstech.}, 45:362--371, 1995.

\bibitem{SchM02}
A.J.~{van der} Schaft and B.M. Maschke.
\newblock {H}amiltonian formulation of distributed-parameter systems with
  boundary energy flow.
\newblock {\em J. Geom. Phys.}, 42:166--194, 2002.

\bibitem{SchM13}
A.J.~{van der} Schaft and B.M. Maschke.
\newblock Port-{H}amiltonian systems on graphs.
\newblock {\em SIAM J. Control Optim.}, 51:906--937, 2013.

\bibitem{SchM23}
A.J. van~der Schaft and V.~Mehrmann.
\newblock Linear port-{H}amiltonian {DAE} systems revisited.
\newblock {\em System Control Let.}, 177:105564, 2023.
\newblock https://doi.org/10.1016/j.sysconle.2023.105564.

\bibitem{Van79}
P.~{Van~Dooren}.
\newblock The computation of {K}ronecker's canonical form of a singular pencil.
\newblock {\em Linear Algebra Appl.}, 27:103--121, 1979.

\bibitem{WeiWS94}
H.~Weiss, Q.~Wang, and J.~L. Speyer.
\newblock System characterization of positive real conditions.
\newblock {\em IEEE Trans. Automatic Control}, 39(3):540--544, 1994.

\bibitem{Wil72b}
J.~C. Willems.
\newblock Dissipative dynamical systems -- {P}art {II}: {L}inear systems with
  quadratic supply rates.
\newblock {\em Arch. Ration. Mech. Anal.}, 45:352--393, 1972.

\end{thebibliography}
\end{document}